\documentclass[11pt]{amsart}

\usepackage{import}

\usepackage{set}

\newcounter{dawidcomments}

\usepackage[T1]{fontenc}

  \title{On the profinite rigidity of free and surface groups}
\author{Ismael Morales}
\date{}

\UseRawInputEncoding

\begin{document}
\maketitle

\begin{abstract}
Let $S$ be either a free group or the fundamental group of a closed hyperbolic surface. We show that if $G$ is a finitely generated residually-$p$ group with the same pro-$p$ completion as $S$, then two-generated subgroups of $G$ are free.  This generalises (and gives a new proof of) the analogous result of Baumslag for parafree groups. Our argument relies on  the following new ingredient: if $G$ is a residually-(torsion-free nilpotent) group and $H\leq G$ is a virtually polycyclic subgroup, then $H$ is nilpotent and the pro-$p$ topology of $G$ induces on $H$ its full pro-$p$ topology. 

Then we study applications to profinite rigidity. Remeslennikov conjectured that  a finitely generated residually finite $G$ with profinite completion $\hat G\cong \hat S$ is necessarily $G\cong S$. We confirm this when $G$ belongs to a class of groups $\Hia$ that has a finite abelian hierarchy starting with finitely generated residually free groups. This strengthens a previous result of Wilton that relies on the hyperbolicity assumption. Lastly, we  prove that the group $S\times \Z^n$ is profinitely rigid within finitely generated residually free groups. 
\end{abstract}

\section{Introduction}
A fundamental idea in the study of infinite groups consists of understanding how much algebraic structure of a group $G$ can be detected in its collection of finite quotients. We address the following question: Suppose that $S$ is either a free group or the fundamental group of a closed hyperbolic surface. Given a finitely generated residually finite group $G$   with the same profinite completion $\hat G\cong \hat S$, what can be said about $G?$  Our  main result is \cref{Titsintro} and provides a partial answer which concerns the subgroup structure of $G$.  Before stating it, we  remark that a recent breakthrough of Jaikin-Zapirain \cite[Theorem 1.1]{And212} implies that such $G$ is residually-$p$, so a group $G$ as above will lie under the assumptions of the following theorem.

\begin{thmA} \label{Titsintro} Let $p$ be a prime and let $S$ be either a free or hyperbolic surface group. Suppose that $G$ is a finitely generated residually-$p$ group with the same pro-$p$ completion as $S$. Then two-generated subgroups of $G$ are free.
\end{thmA}

This is proven in \cref{Tits}. \cref{Titsintro}  generalises (and gives a new proof of) the analogous result of Baumslag \cite[Theorem 4.2]{Bau69} for parafree groups $G$. Recall that a finitely generated group $G$ is termed {\te parafree} (resp. {\te parasurface})  if it is residually nilpotent and its quotients by the terms of its lower central series are the same as those of some free group (resp. surface group). Parafree groups were introduced by Baumslag  \cite{Bau67} and these seem to resemble free groups in many of their structural properties. The proof of Baumslag relies on the fact that a parafree group  can be embedded  into the unit group of a ring of power series $\Z\llbracket X\rrbracket $, an approach which does not seem to carry over to parasurfaces. Our proof of \cref{Titsintro} relies on different cohomological and separability arguments, including our next theorem. 
\begin{thmA} \label{RNsep} Let $G$ be a finitely generated residually-(torsion-free nilpotent) group and let $H\leq G$ be a virtually polycyclic subgroup. Then $H$ is nilpotent and, for all primes $p$, the pro-$p$ topology of $G$ induces on $H$ its full pro-$p$ topology. 
\end{thmA}
This is proven in \cref{pemb}. \cref{RNsep} provides a fairly general setting in which nilpotent subgroups are witnessed by the pro-$p$ completion. This is also crucial ingredient of \cite[Theorem E]{Fru23}, where Fruchter and the author show  that direct products of free and surface groups are profinitely rigid amongst finitely presented residually free groups. We survey during \cref{7.3} more applications that motivate the study of pro-$p$ topologies of infinite groups. The proof of \cref{RNsep} relies on separability properties of nilpotent groups that are particular to this class (at least amongst the class of virtually polycyclic groups), as our next result exhibits. 
\begin{thmA}\label{induced}  Let $G$ be a virtually polycyclic group. Then the following two statements are equivalent: 
\begin{itemize}
    \item The group $G$ is nilpotent. 
    \item For all subgroups $H\leq G$ and all primes $p$, the pro-$p$ topology of $G$ induces on $H$ its full pro-$p$ topology.
\end{itemize}
\end{thmA}

One implication is proven in \cref{pMalcev} and the other in \cref{reci}. Virtually polycyclic groups have been known to be subgroup separable for long 
(after the work of Hirsch and Malcev \cite[Chapter 1]{Seg83}). However, the situation about their pro-$p$ topology does not seem to be as clear and, to  the best of the author's knowledge, this result has not appeared explicitly in the literature before.  

We now put \cref{Titsintro} in more context regarding previous research on relating subgroup structure and profinite completions, which so far has been done specially for three-manifolds. 
For example, a remarkable result is the profinite recognition of the topology of geometric three-manifolds by Wilton--Zalesskii \cite{Wilton17}, where they show that a compact orientable aspherical three-manifold $M$ is hyperbolic if and only if $\hat{\pi_1 M}$ contains no $\Z_p^2$. This was generalised  by Zalesskii \cite[Theorem 1.1]{Zal22}, showing  that a relatively hyperbolic virtually compact special group $G$ is hyperbolic if and only if its profinite completion $\hat G$ contains no $\Z_p^2$. In this sense, one of the difficulties of \cref{Titsintro} is relating subgroup obstructions of an abstract group to subgroup obstructions  of its  pro-$p$ completion in the flavour of the previous results but in the generality of finitely generated residually-$p$ groups, where there is a priori no hierarchy at our disposal. 

Note that \cref{Titsintro} is not immediate even for residually free groups. An elementary property of residually free groups is that non-abelian two-generated subgroups are free (see \cite[Section 4]{Bau62}). However, these can clearly contain the group $\Z^2$ and, even in this context, one needs an additional argument to rule it out. One may expect the analogous result of \cref{Titsintro} to hold for limit groups. As observed in \cref{Titslimit}, our methods  establish  \cref{Titsintro} for the $p$-genus of more hyperbolic limit groups, other than free and surface groups. However, we cannot solve the general case.

\begin{question} Let $p$ be a prime and let $G$ be a group in the $p$-genus of a hyperbolic limit group $L$. Are two-generated  subgroups of $G$ free?
\end{question}

An important property of our argument is that, when $S$ is a free or surface group, $S_{\hp}$ is residually-(torsion-free nilpotent), proven in Propositions \ref{tfnil} and \ref{tfnil2}. For a general limit group $L$, this feature about $L_{\hp}$ is not known. Kochloukova--Zalesskii  study in \cite[Theorem 4.2]{Koc11} the residual properties of a class of pro-$p$ analogues of limit groups (which includes, for example, pro-$p$ completions of ICE groups). However, it is an open problem to determine whether their class contains the pro-$p$ completions of all limit groups (we refer to \cite[Section 9]{Koc11} for more open questions).
 \subsection{Applications to profinite rigidity}
 It is a well-known conjecture of Remeslennikov  that free groups are profinitely rigid (see, for example, \cite{Khu14}*{Question 5.48} and \cite{Nos79}*{Question 12}). This can be naturally restated to include surface groups as follows.
\begin{conjecture} \label{Rem} Let $S$ be either a free or a surface group and let $G$ be a finitely generated residually finite group with $\hat G\cong \hat S$. Then  $G\cong S$.
\end{conjecture}
 
Wilton \cites{Wil18, Wil21} solved \cref{Rem} for limit groups $G$. 
There is a second proof of \cref{Rem} when $S$ is a surface group and $G$ is a finitely generated residually free group by Fruchter and the author \cite[Corollary D]{Fru23}. In \cref{Remthm}, we establish \cref{Rem} when $G$ belongs to the following abelian hierarchy, which naturally includes limit groups (since these admit an abelian hierarchy terminating in free groups \cite{Sel01}*{Theorem 4.1}).

\begin{defi} \label{defhia} Let $\Hia$ denote the smallest class of groups such that 
\begin{enumerate}
    \item $\Hia$ contains all finitely generated residually free groups. 
    \item If $A, B$ are in  $\Hia$ and $C$ is a virtually abelian subgroup of $A$ and $B$, then the amalgamated product $A *_C B$ belongs to $\Hia$. 
    \item If $A$ is in  $\Hia$ and $\theta\colon C\lrar C'$ is an isomorphism between  virtually abelian subgroups of $A$, then the HNN extension $A\,  *_{C, \theta}$ belongs to $\Hia$.
    \item If $A$ is in  $\Hia$ and $A$ embeds as a finite-index subgroup in $B$, then $B$ belongs to $\Hia$.
\end{enumerate}
\end{defi}

Once we know that a group $G$ as in \cref{Rem} does not have Baumslag--Solitar subgroups by \cref{Titsintro} and \cite[Theorem 1.1]{And212}, $G$ will be hyperbolic by the Bestvina--Feighn combination theorem \cite{Bes92}. From this, we follow the outline of the proof of Wilton (while also incorporating  Wise--Haglund's separability and structural properties of    hyperbolic groups with a malnormal
quasi-convex hierarchy \cites{Hag08, Wis21}) to prove the following generalisation. 
\begin{corA}\label{Remthm} \cref{Rem} holds for groups $G$ belonging to the class $\Hia$.
\end{corA}
This is proven in Theorems \ref{Remfree} and \ref{Remsur}. There are several reasons why it is natural to study $\Hia$, other than being a family that naturally incorporates other groups for which \cref{Rem} was solved. 
The consideration of hierarchies with abelian edges has provided a fruitful setting in which  combination theorems for relevant classes of groups have been established:   for coherent groups \cite[Corollary 5.4]{Thu20}, for hyperbolic groups \cite{Bes92}, for free groups \cites{She55, Swa86, Sta91}  and for parafree groups \cite{And21}. We pose a question that already seems to contain most of the difficulties that \cref{Rem} presents for more complicated group hierarchies.

\begin{question} Let $F$ be a finitely generated free group and let $H$ be a finitely generated malnormal subgroup. Consider the double $\Ga=F*_{H=\ove H} \ove F$ of the group $F$ along its subgroup $H$. In particular, $\Ga$ is hyperbolic and virtually compact special by \cite[Corollary 5.3]{Git97} and \cite[Corollary B]{Hsu15}.  Is \cref{Rem}  true for $\Ga$?
\end{question}

Another class of groups where the conclusion of \cref{Titsintro} ensures hyperbolicity is one-relator groups (by the advances of Louder--Wilton \cites{Lou22, Lou22b} and Linton \cite{Lin22}). \cref{Remthm} covers some one-relator groups such as amalgamated products or HNN extensions of free groups along cyclic subgroups, which have the form $\lan F_1, F_2\, |\, u_1=u_2\ran$ and $\lan F, t \, |\,tut^{-1}=v\ran$. Nevertheless, we do not know whether residually finite one-relator groups $G$ with the same profinite completion as a free or surface group $S$ necessarily satisfy $G\cong S$. Opportunely, \cref{hvspecialintro} at least ensures that such $\Ga$ are  hyperbolic and virtually special. 

\begin{corA} \label{hvspecialintro} Let $p$ be a prime and let $G$ be a one-relator group in  $p$-genus of a free or hyperbolic surface group. Then $G$ is hyperbolic and virtually special. 
\end{corA}

 Again we stress that \cref{hvspecialintro} extends to other genera by \cite[Theorem 1.1]{And212}.  Our last result is \cref{prods}, which provides infinitely many examples of groups that are determined within the class of residually free groups by their finite quotients. We extend this class of groups using similar methods in \cite{Mor23} by looking at examples coming from three-manifold groups. 

\begin{thmA}\label{prods}  Let $G$ be a finitely generated residually free group and let $S$ be a free or surface group. Suppose that $\hat G\cong \hat S\times \hat \Z^n$. Then $G\cong S\times \Z^n$. 
\end{thmA}

We prove this theorem in \cref{prodsec}. The first step consists on showing that $G$ and $S\times \Z^n$ have isomorphic centre, for which we use the theory of $L^2$-Betti numbers. This does not occur in general because  Lubotzky proved in \cite[Proposition 1.5]{Lub14} that having a non-trivial centre is not a profinite property. As we mentioned above, our methods give an extension of \cref{prods} to direct products of free and surface groups, under the extra assumption of $G$ being finitely presented, in a joint work with Fruchter \cite[Theorem E]{Fru23}.

\subsection{Acknowledgements} The author is funded by the Oxford--Cocker Graduate Scholarship. This work has also received funding from the European Research Council (ERC) under the European Union's Horizon 2020 research and innovation programme (Grant agreement No. 850930). I am grateful to Andrei Jaikin-Zapirain introducing me to these questions on profinite rigidity and for suggesting \cref{induced}. I would like to thank Jonathan Fruchter, for many insightful conversations; Sam Hughes and Gareth Wilkes, for multiple  corrections and suggestions; and also Alejandra Garrido and Pavel Zalesskii, for a helpful correspondence. Lastly, special thanks are due to my supervisors Dawid Kielak and Ric Wade for their  stimulating guidance and for carefully reading this paper.

\subsection{Organisation of the paper} 
In  \Cref{proCsec} we revise some properties of profinite (and pro-$p$) groups and in \Cref{nilpotent} we study the induced pro-$p$ topology on subgroups of nilpotent groups, proving Theorems \ref{RNsep} and \ref{induced}. Then, we show \cref{Titsintro} in \cref{Titssec}. Next, we review $L^2$-Betti numbers in \Cref{Betti} and the structure of residually free groups in \cref{limitsec} to prove \cref{prods} in \cref{prodsec}. In subsequent sections, we review the tools that are needed to establish \cref{Remthm} and \cref{hvspecialintro}. 

\section{Pro-\texorpdfstring{$\CC$}{C} groups, topologies and completions} \label{proCsec}

A non-empty class of finite groups $\CC$ is a {\te formation} if it is closed under taking homomorphic images and subdirect products. 
Given a finitely generated group $ \La$, we say that it is {\te residually-$\CC$} if, for all $1\neq \lambda \in \La$, there exists $N\n \La$ such that $\lambda \notin N$ and $\La/N \in \CC$. In this case, following the terminology of \cite{Gru11}, we define the  {\te $\CC$-genus} of $\La $ as the collection of finitely generated residually-$\CC$ groups $\Ga$ with  isomorphic pro-$\CC$ completion $\Ga_{\hat\CC}\cong  \La_{\hat\CC}$. Equivalently, the latter happens if $\Ga$ and $\La$ have the same collection of isomorphism types of finite quotients belonging to $\CC$ (\cref{proC}).  We will mostly think of $\CC$ as being the formations of all finite groups,  all finite soluble groups,  all finite nilpotent groups or all  finite $p$-groups (for some prime $p$). In these cases, we refer to the $\CC$-genus as the finite,  soluble,  nilpotent and $p$-genus, respectively. Lastly, we say that a finitely generated residually finite $\La$ is {\te profinitely rigid} if the only $\Ga$ in its finite genus is $\Ga\cong \La$.  


\begin{notation} Given an abstract group $\Ga$, we denote its profinite completion (resp. pro-$p$ completion) by $\hat \Ga$ (resp. $\Ga_{\hp}$). Bold letters $\GG, \HH, \KK, \UU$ will denote profinite groups. In addition, we write $\HH\lc \GG$ or $\HH\lo \GG$ (resp. $\HH \nc \GG$ or $\HH \no \GG$) to indicate that the subgroup $\HH\leq \GG$ (resp. normal subgroup $\HH\n \GG$) is closed or open.  
\end{notation} 
It will be important to keep in mind during the rest of the exposition the following theorem, which has already been discussed during the introduction.
\begin{thm}[Theorem 1.1, \cite{And212}] \label{And} Let $S$ be a finitely generated free or surface group. Suppose that $G$ belongs to the finite, soluble or nilpotent genus of $S$. Then $G$ belongs to the $p$-genus of $S$ for all primes $p$.
\end{thm}
These inclusions between genera are not true in general, since Lubotzky \cite{Lub14} showed that the properties of being residually soluble or residually-$p$ are not  profinite invariants.

\subsection{Pro-$\CC$ completions} We now fix a formation $\CC$ of finite groups. A pro-$\CC$ group can be defined as a profinite group $\GG$ such that for all its open normal subgroups $\UU\no \GG$, one has that $\GG/\UU\in \CC$. Alternatively, a pro-$\CC$ group is just an inverse limit of groups belonging to $\CC$ that is endowed with the profinite topology. 

In order to relate residual properties of abstract groups and profinite invariants, we recall the pro-$\CC$ completion functor. Given an abstract group $G$, its collection of quotients belonging to $\CC$ is encoded in its  pro-$\CC$ completion $G_{\hat \CC}$. As an abstract group, it is defined  to be the inverse limit of finite groups
\[G_{\hat \CC}=\varprojlim_U G/U\]
where $U$ ranges over finite-index normal subgroups of $G$ such that $G/U\in \CC$. We always consider $G_{\hat \CC}$ endowed with the profinite topology, which turns it into a pro-$\CC$ group. There is a natural and canonical group homomorphism $\iota_{\hat \CC}\po G\lrar G_{\hat \CC}$ and, in these terms, a group $G$ is residually-$\CC$ if and only if $\iota_{\hat \CC}$ is injective. This map $\iota_{\CC}$ enjoys the following universal property (which, in fact, provides a more functorial definition of $G_{\hat \CC}$).

\begin{prop} \label{uniprop} Let $\Ga$ be an abstract group. Let $\iota_{\hat{\CC}}\po \Ga\lrar \Ga_{\hat{\CC}}$ be the canonical map. The data $(\Ga_{\hat  \CC}, \iota_{\hat \CC})$ is characterised by the following properties. 
\begin{enumerate}
    \item $\Ga_{\hat{\CC}}$ is a pro-$\CC$ group. 
    \item The map $\iota_{\hat{\CC}}$ has dense image. 
    \item For any pro-$\CC$ group $\KK$ and every group homomorphism $f\po \Ga\lrar \KK$ with dense image, there exists a continuous homomorphism $f_{\hat{\CC}}\po \Ga_{\hat{\CC}}\lrar \KK$ such that the diagram 
    \begin{equation*}
        \begin{tikzcd}
         \Ga \ar[dr, "f"] \ar[r, "\iota_{\hat \CC}"] & \Ga_{\hat{\CC}}\ar[d, dashrightarrow, "f_{\hat{\CC}}"]\\
           & \KK
        \end{tikzcd}
    \end{equation*}
    is commutative. 
\end{enumerate}
When we say that the pair $(\Ga_{\hat \CC}, \iota_{\hat \CC})$ is unique, we mean the following: if  $\HH$ is a pro-$\CC$ group and $\iota\po \Ga\lrar \HH$ is a group homomorphism such that the pair $(\HH, \iota)$ verifies the three properties above, then there exists an isomorphism of pro-$\CC$ groups $\al\po \Ga_{\hat\CC}\lrar \HH$ such that $\al\circ \iota_{\hat \CC}=\iota$. 
\end{prop}

We stress the fact that the previous discussion does not rely on whether $\iota_{\hat \CC}$ is injective (equivalently $G$ is residually-$\CC$) or not. A priori the pro-$\CC$ completion only captures information about finite quotients of $G$ in $\CC$ and not about how we can distinguish elements of $G$ in these quotients. There are, for example, finitely generated  residually finite groups with free pro-$p$ completion for all primes $p$ that are not residually-$p$ for any prime $p$, such as the trefoil knot group $\lan a, b\, | a^2b^3\ran$. The philosophy that finite quotients of a group are encoded in their profinite completion is  precised in the following folkloric theorem. 

\begin{thm}[Dixon--Formanek--Poland--Ribes \cite{Dix82}]\label{proC} Let $\CC$ be a formation of finite groups. Two abstract groups $\Ga$ and $\Lambda$ have the same class of isomorphism types of quotients belonging to $\CC$ if and only if $\Ga_{\hat{\CC}}\cong \Lambda_{\hat{\CC}}$. 
\end{thm}

So we can completely reformulate in terms of the pro-$\CC$ completions what it means for a group to be in the $\CC$-genus of another one. This is important because now we have turned a bare collection of finite groups into a group and now we can do, say, group cohomology with this information. Crucially, this functor ``pro-$\CC$ completion'' is right-exact.

\begin{prop}[\cite{Rib00}, Proposition 3.2.5] \label{secC} Let $\CC$ be a formation that is closed under taking normal subgroups. Given a exact sequence of abstract groups of the form $N\lrar G\lrar Q\lrar 1$, the induced sequence of maps 
$N_{\hat \CC}\lrar G_{\hat \CC}\lrar Q_{\hat \CC}\lrar 1$ is exact. 
\end{prop}

 Given an injective map $N\lrar G$, the question of whether the induced map $N_{\hat \CC}\lrar G_{\hat \CC}$ is also injective will be a crucial topic of discussion later in this section and during \cref{nilpotent}, especially with the aim of proving \cref{Titsintro}. On the other hand, \cref{secC} also tells us that pro-$\CC$ completions are better understood when studying quotients than when studying subgroups, as the following direct consequence exhibits. 

\begin{prop}[{\te Pro-$\CC$ completion of a quotient}] \label{propquo} Let $G$ be an abstract group and let $N\n G$ be a normal subgroup. Suppose that $\CC$ is a formation as in \cref{secC}. Denote by $\iota\po N\lrar G$ the inclusion. Then the canonical and natural map
\[\left(G/N\right)_{\hat \CC} \lrar G_{\hat \CC}\Big/ \ove{\iota_{\hat \CC}(N)} \]
is an isomorphism. 
\end{prop}

\subsection{The notion of index in profinite groups}

 Given a profinite group $\GG$ and a subgroup $\HH$, we will not need for our purposes to define the formal meaning of $|\GG: \HH|$ as a {\te supernatural number} (as in, for example, \cite[Section 2.3]{Rib00}).

However, recall we say that $p^{\infty}$ divides
the index $|\GG: \HH|$ if and only if for all $n$ there exists a normal open $\UU\no \GG$ such that $p^n$ divides the index of $\HH \UU/\UU$ in the finite group $\GG/\UU$. Observe that, if $\GG$ is a pro-$p$ group, then $ \HH\lc \GG$ has infinite index if and only if $p^{\infty}$ divides the index $|\GG: \HH|$. However, for general profinite groups, the latter is much stronger, as the following example shows.

\begin{eg}\label{Zhat} The profinite group $\hat \Z$ is canonically isomorphic to $\prod_p \Z_p$, where $\Z_p$ denotes the additive part of the ring of $p$-adic integers (isomorphic to the pro-$p$ completion of $\Z$). So $\hat \Z$ surjects in a natural way onto the profinite group $\prod_p \Z/p\Z$, where the kernel $\KK$ is isomorphic to $\hat \Z$ and will have an index in $\hat \Z$ equal to the supernatural number $\prod_p p$. In particular, no prime $p$ will satisfy that $p^2$ divides the index $|\hat \Z: \KK|$. 
\end{eg}

\subsection{Nilpotent groups and lower central series} \label{lower} The terms of the {\te lower central series} of a group $G$ are defined recursively as $\gamma_1 G=G$ and $\gamma_{n+1} G=[\gamma_n G, G]$ if $n\geq 1$. Similarly, the terms of the $p$-{\te lower central series} of $G$ are defined as $\gamma_{1,p}G=G$ and $\ga_{n+1,p}G=[ \ga_{n,p}G, G]\ga_{n,p} G^p$ if $n\geq 1$. The expression $\ga_n \Ga_{\hp}$ may be ambiguous. We clarify that it means $\ga_n G$ with $G=\Ga_{\hp}$,  defined as a subgroup of $\Ga_{\hp}$ (so we do not mean to denote the pro-$p$ completion of $ \ga_n \Ga$). 

 A group $G$ is {\te nilpotent} if there exists an integer $c\geq 0$ such that $\ga_{c+1}G=1$.
 For example, abelian groups are nilpotent groups of nilpotency class 1. We will study the profinite and pro-$p$ topology of nilpotent groups in Sections \ref{inducedtop} and \ref{nilpotent}.

\begin{prop}[\cite{Dix99}, Proposition 1.19 and Exercise 17] \label{gak}
Let $\GG$ be a finitely generated pro-$p$ group. Then $\ga_k  \GG$ is a closed subgroup of $\GG$. 
\end{prop} 
In fact, the previous is true for finitely generated profinite groups, although this is a much deeper result due to Nikolov and Segal \cite[Theorem 1.4]{Nik07}.

\begin{lemma} \label{gakp} Let $\Ga$ be a finitely generated group and let $k\geq 0$ be an integer. Then the canonical map $\Ga/\ga_k \Ga\lrar \Ga_{\hp}/\ga_k\Ga_{\hp}$ is a pro-$p$ completion map. 
\end{lemma}
\begin{proof}
The group $\ga_k\Ga_{\hp}$ are closed subgroups of $\Ga_{\hp}$ by \cref{gak} and  contain $\ga_k\Ga$ as a dense subgroup, so $\ga_k\Ga_{\hp}$ is the closure of $\ga_k\Ga$. By \cref{propquo}, the conclusion follows. 
\end{proof}

\subsection{Induced profinite topologies} \label{inducedtop} Given a formation $\CC$ and a group $G$ the {\te full pro-$\CC$ topology} of $G$ is defined to be the initial topology of the map $G\lrar G_{\hat \CC}$. In other words, it is the smallest topology that contains among its open subsets all the cosets of $U\n G$ such that $G/U\in \CC$. When $\CC$ is the formation of all finite groups (resp. finite $p$-groups), we call it the {\te profinite (resp. pro-$p$) topology}. 

\begin{notation} When the formation $\CC$ is understood from the context, we write $H\lc G$ or $H\lo G$ (resp. $H\nc G$ or $H\no G$) to denote that the subgroup $H\leq G$ (resp. normal subgroup $H\n G$) is closed or open in the pro-$\CC$ topology.  
\end{notation}

A useful result that  relates the open subgroups of a group $G$ with the open subgroups of $\hat G$  in a  explicit way is the following correspondence.

\begin{prop}[\cite{Rib00}, Proposition 3.2.2] \label{bij} Let $G$ be a residually-$\CC$ group and identify $G$ with its image in $G_{\hat \CC}$. 
\begin{itemize}
    \item[(a)] Let \[\The: \{H \, |\, H\lo G\}\lrar \{\HH\, |\, \HH\lo G_{\hat \CC}\}\] be the map that assigns to each open subgroup $H\lo G$ its closure $\ove{H}$ in $ G_{\hat \CC} $. Then $\The$ is a bijective correspondence with inverse given by $\The^{-1}(\HH)=\HH\cap G$. 
    \item[(b)] Given $H\lo G$,  it is true that $H\no G$ if and only if $\The(H)\no G_{\hat \CC}$.
    \item[(c)] Given $H\lo G$, the index $|G: H|$ is equal to the index $|G_{\hat \CC}: \The(H)|$.
    \item[(d)] The canonical map $H\lrar \The(H)$ is a pro-$\CC$ completion map. 
\end{itemize}
\end{prop}

\begin{defi} \label{embdefi} Let $\CC$ be a formation of finite groups. Assume that $H\leq G$. Denote by $\iota\po H\lrar G$ the injection. We say that $H$ is {\te pro-$\CC$ embeddable into} $G$ if the induced map 
\[\iota_{\hat \CC}\po H_{\hat \CC}\lrar G_{\hat \CC} \]
is an injection. 
\end{defi}
Observe that this is a transitive relation in the sense that whenever $A\leq B\leq C$, with $A$ being pro-$\CC$ embeddable into $B$ and $B$ being pro-$\CC$ embeddable into $C$, then $A$ is pro-$\CC$ embeddable into $C$. Other than stating \cref{bij},  we introduced the pro-$\CC$ topology on abstract groups because it allows us to reformulate pro-$\CC$ embeddability intrinsically in $G$ as follows. 

\begin{lemma}[See \cite{Rib00}, Lemma 3.2.6] \label{embed} Let $H\leq G$ be two groups. Then $H$ is pro-$\CC$ embeddable into $G$ if and only if the pro-$\CC$ topology of $G$ induces on $H$ its full pro-$\CC$ topology. 
\end{lemma}

\begin{obs}\label{strat} If $H$ is finitely generated, in order to prove that the pro-$\CC$-topology of $G$ induces the full pro-$\CC$ topology of $H$, it suffices to check the following: For any characteristic subgroup $\La\n H$ with $H/\La\in \CC$, there exists $\Ga\n G$ with $G/\Ga\in \CC$ such that $\Ga\cap H\sub \La$.
\end{obs}

We will show in \cref{pMalcev} that subgroups of finitely generated nilpotent groups are pro-$p$ embeddable. We finish this section by analysing some of its consequences  (for instance, the fact that Hirsch lengths will be pro-$p$ invariants of finitely generated torsion-free nilpotent groups). 

\begin{prop}\label{tfnil} Let $G$ be a finitely generated group such that the central quotients $G/\ga_k G$ are torsion-free for all $k$. Then, for all primes $p$, $G_{\hp}$ is residually-(torsion-free nilpotent). In fact, the quotients $G_{\hp}/\ga_n G_{\hp}$ are torsion-free for all $n$. 
\end{prop}

\begin{proof} 
Since $G_{\hp}$ is residually-$p$,  $\bigcap_k \ga_{k, p} G_{\hp}=\{1\}$ and, in particular, $\bigcap_k \ga_{k}G_{\hp} =1$. It remains to prove that $G_{\hp}/\ga_k G_{\hp}$ is torsion-free. 
 We know that $\Ga=G/\ga_{c+1} G$ is  a finitely generated torsion-free nilpotent group. We are going to prove that this is enough to conclude that its pro-$p$ completion $\Ga_{\hp}$ is torsion-free. The proof of the proposition would then be complete because the canonical map $\Ga_{\hp}\lrar G_{\hp}/\ga_k G_{\hp}$ is an isomorphism by \cref{propquo}.

Now let $\Ga$ be a  finitely generated torsion-free nilpotent group. Then $\Ga/Z(\Ga)$ is also torsion-free (see \cite[Corollary 2.22]{Cle17}). Then, iterating the previous observation, we deduce that the lower central series \[1=\ga_{c+1}\Ga\unlhd \ga_c\Ga\n \cdots \ga_2\Ga\unlhd \ga_1\Ga=\Ga\] of $\Ga$ will be such that $\ga_{k}\Ga/\ga_{k+1}\Ga\cong \Z^{n_k}$ for all $1\leq k\leq c$ and some non-negative integers $n_1, \dots, n_c$.

 Notice that the canonical maps $(\ga_k \Ga)_{\hp}\lrar \ga_k \Ga_{\hp}$ are isomorphisms (for example, by \cref{pMalcev}). From this and \cref{propquo}, we derive that  $(\ga_{k+1}\Ga/\ga_k \Ga)_{\hp}\lrar \ga_{k+1}\Ga_{\hp}/\ga_k\Ga_{\hp}$ are isomorphisms of pro-$p$ groups for all $1\leq k\leq c$ (recall that these groups are, by \cref{gakp}, the pro-$p$ completions of $\ga_{k+1}\Ga/\ga_k \Ga$)
 So there is a lower central series of $\Ga_{\hp}$ by closed subgroups
\[1=\ga_{c+1}\Ga_{\hp}\unlhd \ga_c\Ga_{\hp}\unlhd \cdots \unlhd \ga_2\Ga_{\hp}\unlhd \ga_1\Ga_{\hp}=\Ga_{\hp},\]
for which each quotient $\ga_{k+1}\Ga_{\hp}/\ga_k \Ga_{\hp}\cong \Z_p^{n_k}$ is torsion-free. So $\Ga_{\hp}$ is torsion-free, as we claimed before.
\end{proof}

It is not hard to see that the class of groups satisfying the assumptions of \cref{tfnil} is closed under taking free products. However, the following proposition ensures many one-ended examples as well. 

\begin{prop} \label{tfnil2} Let $G$ be isomorphic to a RAAG or to a one-relator group $F/\lan\!\lan r\ran\!\ran$ such that if $k$ is the biggest non-negative integer such that $r\in \ga_k F$, then $r\cdot \ga_{k+1}F$ is not a proper power in $\ga_{k}F/\ga_{k+1}F$. In these cases, $G$ satisfies the assumptions of \cref{tfnil}. In particular, this is true if $G$ is a  surface group. 
\end{prop}
\begin{proof} When $G$ is a RAAG, this is   \cite[Theorem 6.4]{Wade15}. When $G=F/\lan\lan r\ran\ran$, this proposition is a direct consequence of \cite[Section 1]{Lab70}. The reason is the following. The main theorem of \cite{Lab70} states that, under our assumptions on the relator $r$, the  graded algebra  $\bigoplus_n \ga_n G/\ga_{n+1} G$ is a free $\Z$-Lie algebra. In particular, the subsequent quotients $\ga_n G/\ga_{n+1} G$ are free abelian  (the author also computes their rank). Since surface groups are defined by a single relator that satisfies the assumptions of $r$ in \cref{tfnil2}, the last conclusion follows. 
\end{proof}
We observe that the class of one-relator groups of \cref{tfnil2} includes many limit groups. For example, those arising as doubles $F*_{u=\bar{u}} \bar{F}$ of a free group $F$ along a relator ${u}\in F$ that satisfies the assumptions of $r$ in \cref{tfnil2}. We remark in \cref{Titslimit} how this can be used to establish \cref{Titsintro} for the $p$-genus of many  hyperbolic limit groups, other than free and non-abelian surface groups. 

\section{Profinite and pro-\texorpdfstring{$p$}{p} topologies on polycyclic groups} \label{nilpotent}
We return to analysing the structure and pro-$p$ topology of nilpotent groups, as we initiated in \cref{lower}. We prove \cref{pMalcev} and \cref{reci} (which constitute \cref{induced} from the introduction) and then show \cref{pemb} (\cref{RNsep}), which allows to witness polycyclic subgroups from the pro-$p$ completion and plays a crucial role in the proof of \cref{Titsintro} during \cref{Titssec}.

It was proven by Hirsch \cite[Theorem 3.25]{Hir46} that polycyclic groups are residually finite, and Malcev later in 1948 proved that they are, in fact, LERF \cite[Exercise 11, Chapter 1]{Seg83}. The former result of Hirsch has a pro-$p$ variant for nilpotent groups, given by Gruenberg (\cref{resp}), while Malcev's result does not. Still,  it follows from Malcev's theorem the weaker property that the induced topology on subgroups of polycyclic groups is the full profinite topology. We study in \cref{pMalcev} a pro-$p$ version of this statement. Before stating it, we first review Gruenberg's result and some definitions. 

Given a set of primes $P$ and a non-negative integer $n$, we say that $n$ is a $P$-{\te number} if all the primes dividing $n$ are contained in $P$. We say that a group $G$ is {\te residually}-$P$ if, for all $1\neq g\in G$, there exists a prime $p\in P$ and a finite-index normal subgroup $U\n G$, where $|G: U|$ is a power of $p$, such that $g\notin U$. Notice that, by Cauchy's theorem, the set of primes $P$ appearing as the orders of elements of a group $G$ is exactly  the smallest set of primes $P$ such that every element with finite order in $G$ has a $P$-number order. 

\begin{thm}[\cite{Gru57}, Theorem 2.1] \label{resp} Let $G$ be a finitely generated nilpotent group.
\begin{itemize}
    \item[(a)] If $G$ is torsion-free, then $G$ is residually-$p$ for all primes $p$. 
    \item[(b)] Otherwise, let $P$ be the (non-empty) collection of primes $p$ that appear as orders of elements in $G$.  Then $G$ is residually-$P$. 
\end{itemize}
\end{thm}

\subsection{Induced pro-$p$ topology on nilpotent groups}  Continuing the discussion at the beginning of the section, one may also try to extend the result of Malcev to the pro-$p$ setting by, in the first place, figuring out a definition of $p$-LERF groups. A tentative definition is to say that a group $G$ is {\te $p$-LERF} if all its finitely generated subgroups are closed in the pro-$p$ topology. However, this may not be an interesting notion for infinite groups, since a virtually polycyclic group is $p$-LERF if and only if it is a finite $p$-group (see \cref{plerf}). In any case, many infinite groups have the weaker property that the induced pro-$p$ topology on subgroups is their full pro-$p$ topology. We confirm this in the next theorem.

\begin{thm}\label{pMalcev} Let $G$ be a finitely generated nilpotent group and let $H$ be any subgroup. Then the pro-$p$ topology of $G$ induces on $H$ its full pro-$p$ topology.
\end{thm}

By \cref{embed}, the previous condition on each $H\leq G$ is equivalent to the injectivity of the induced map $H_{\hp}\lrar G_{\hp}$. Recall in this case we say that $H$ is pro-$p$ embeddable into $G$ (\cref{embdefi}).  \cref{pMalcev}  says that the functor ``pro-$p$ completion'' is exact in the category of finitely generated nilpotent groups (since it is always right-exact by \cref{secC}). Before the proof of  \cref{pMalcev}, we will review how one can extract roots of subgroups in nilpotent groups in \cref{rad}.

\begin{defi} Let $G$ be a group, let $P$ be a collection of primes and assume $H\leq G$. We define the $P$-{\te radical} of $H$ in $G$ as the set 
\[R_P(H)=\{g\in G\, \, |  \mbox{ there exists a $P$-number $n$ with $g^n\in H$}\}.\]
\end{defi}

\begin{prop}\label{rad}  Let $G$ be a finitely generated nilpotent group and let $H\leq G$. Then $R_P(H)$ is a subgroup of $G$ and the finite index $|R_P(H): H|$  is a $P$-number. 
\end{prop}
\begin{proof} Let $K$ be the subgroup generated by $R_P(H)$. Since $G$ is noetherian, then $K$ is finitely generated. So $K$ is generated by a finite subset of $R_P(H)$. Then, by \cite[Theorem 2.24]{Cle17}, it follows that $K=R_P(H)$ and that $H$ has a finite index in $R_P(H)$ which is a $P$-number. 
\end{proof}

\begin{proof}[Proof of \cref{pMalcev}] We divide the proof in  two steps. First, we prove in \cref{1} that the conclusion follows if we additionally suppose that $H$ is normal, and then we prove the general case when $H$ is any subgroup. 
\begin{claim} \label{1} The conclusion follows if $H$ is normal in $G$. 
\end{claim}
\begin{proof}[Proof of \cref{1}] We follow the strategy depicted in \cref{strat}. Let $\La\n H$ be a characteristic subgroup of $H$ whose index $|H: \La|$ is finite and a $p$-power. In particular, $\La$ is normal in $G$.  Let $P$ be the collection of primes different from $p$. The subgroup $N=R_P(\La)$ is  normal in $G$ and it contains $\La$ with finite-index. It is immediate to see that $N\cap H=\La$ and that the finite group  $G/N$ does not contain elements of order coprime with $p$. The latter implies, by \cref{resp}, that the nilpotent group $G/N$ is residually-$p$.   We have an injection $HN/N\inc G/N$, where $HN/N$ is a finite subgroup. There exists a finite $p$-power index normal subgroup $K\n G$ containing $N$ such that $K/N$ intersects $HN/N$ trivially (or, equivalently, $K\cap HN=N$). Lastly, we check that $K\cap H=\La$, which would complete the proof of the claim. For this, observe that $K\cap H= (K\cap HN)\cap H=N\cap H=\La. $
\end{proof}

We now finish the proof of the theorem. Consider any subgroup $H\leq G$. Let $c$ denote the nilpotency class of $G$. By \cite[Theorem 10.3.3]{Hal76}, there is an ascending chain $H=H_0\n H_1\n \cdots \n H_c=G$.  By \cref{1}, each $H_i$ is pro-$p$ embeddable into $H_{i+1}$. Finally, by transitivity, $H$ is pro-$p$ embeddable into $G$. 
\end{proof}

\begin{remark} Gareth Wilkes pointed out to the author another proof of \cref{pMalcev}.  The argument can be sketched as follows. Let $H$ be a subgroup of a finitely generated nilpotent group $G$. Since pro-nilpotent groups are direct products of their $p$-Sylows (see \cite[Proposition 2.3.8]{Rib00}), then there is a commutative diagram  
  \begin{equation} 
     \begin{tikzcd}
       H \ar[r, hookrightarrow, "\iota"] \ar[d] & G \ar[d]\\
       H_{\hp} \ar[d, hookrightarrow] \ar[r, "\iota_{\hp}"]& G_{\hp} \ar[d, hookrightarrow] \\
       \hat H \ar[r, hookrightarrow, "\hat \iota"] & \hat G
     \end{tikzcd}
 \end{equation}
where we denoted by $\iota\po H\lrar G$ the injection map and by $\iota_{\hp}$ and $\hat \iota$ the natural induced maps between profinite and pro-$p$ completions, respectively. Since  $G$ is LERF, the map $\hat \iota$ is injective and, consequently, $\iota_{\hp}$ is also injective, as we wanted to prove. 
\end{remark}

\subsection{Induced pro-$p$ topology on polycyclic groups}  \label{7.2}
As we anticipated, although the profinite version of \cref{pMalcev} holds for virtually polycyclic groups, we state the following converse to the pro-$p$ version that shows the limitations of this phenomenon outside of the world of nilpotent groups.

\begin{thm}\label{reci}  Let $G$ be a virtually polycyclic group. Suppose that for every prime $p$ and every subgroup $H\leq G$, the pro-$p$ topology of $G$ induces on $H$ its full pro-$p$ topology. Then $G$ is nilpotent. 
\end{thm}

\cref{reci} in conjunction with \cref{pMalcev} constitute \cref{induced} from the introduction.  Before the proof of \cref{reci}, we start with a chain of preliminary lemmas. 

\begin{lem} \label{rec0} Let $P$ be a finite $p$-group and let $Q$ be a finite $q$-group, with $p$ and $q$ being different primes. Let $\theta \po Q\lrar \Aut(P)$ be a group homomorphism and let $G=P\rtimes Q$ be the corresponding semi-direct product. Then the following are equivalent:
\begin{enumerate}
    \item The map $\theta$ is trivial. 
    \item $G$ is isomorphic to $P\times Q$. 
    \item $G$ is nilpotent.
    \item The subgroup $P$ is pro-$p$ embeddable into $G$. 
\end{enumerate}
\end{lem}
\begin{proof} The equivalence between (1) and (2) is well-known and follows from the definition of semi-direct product (since $Q$ is normal if and only if it is central). By \cite[Theorem 10.3.4]{Hal76}, a finite group is nilpotent if and only it is a direct product of finite $t$-groups with $t$ ranging over prime numbers. So the equivalence between (2) and (3) is also clear. We notice that if (1) is given, then (4) is obviously ensured. So it suffices to prove that (4) implies (1). We consider the trivial subgroup $\{1\}$, which has $p$-power index in $P$ and it is normal. By assumption, there must exist a normal subgroup $N\n G$ of $p$-power index such that $N\cap P\sub \{1\}$, so $N\cap P=\{1\}$. Since the index is coprime to $q$, then $N$ is a $q$-Sylow subgroup of $G$, as well as $Q$. So $N$ and $Q$ are conjugated to each other by the second Sylow's theorem, which shows that $Q$ is also normal and (1) follows. 
\end{proof}

\begin{lem} \label{rec1} Let $p$ be a prime. Suppose that $N \n G$ and $N\leq H\leq G$ have the property that the induced maps $N_{\hp}\lrar H_{\hp}$ and $N_{\hp}\lrar G_{\hp}$ are injective. Then,  the induced map $H_{\hp}\lrar G_{\hp}$ is injective if and only if the induced map $\left(H/N\right)_{\hp}\lrar \left(G/N\right)_{\hp}$ is injective. 
\end{lem}
\begin{proof}  If the induced maps $N_{\hp}\lrar H_{\hp}$ and $N_{\hp}\lrar G_{\hp}$ are injective, we have, by \cref{secC}, a commutative diagram of the form
\begin{equation}
    \begin{tikzcd}
    1\ar[r] & N_{\hp} \ar[r]\ar[d] & H_{\hp}\ar[r]\ar[d] & \left(H/N\right)_{\hp}\ar[r] \ar[d] & 1\\
    1\ar[r] & N_{\hp} \ar[r] & G_{\hp} \ar[r] & \left(G/N\right)_{\hp}\ar[r] & 1,
    \end{tikzcd}
\end{equation}
where all the arrows are the natural ones and the rows are exact. A simple diagram chase gives the desired conclusion. 
\end{proof}

\begin{lem} \label{rec01} Let $P$ be a direct product $P_1\times \cdots\times P_m$, where each $P_i$ is a finite $p_i$-group and $p_i\neq p_j$ if $i\neq j$. Suppose that $q$ is a prime not belonging to the set $\{p_1, \dots, p_n\}$ and let $Q$ be a finite $q$-group. Let $\theta\po Q\lrar \Aut(P)$ be a group homomorphism and let $G=P\rtimes Q$ be the corresponding semi-direct product. Then the following are equivalent:
\begin{enumerate}
    \item The map $\theta$ is trivial. 
    \item $G$ is isomorphic to $P\times Q$. 
    \item $G$ is nilpotent.
    \item The subgroup $P$ is pro-$p_i$ embeddable into $G$ for all $i$. 
\end{enumerate}
\end{lem}
\begin{proof} The equivalence between the first three items is established in the exact same way as in \cref{rec0}. Again, it is clear that $(2)$ implies (4). We suppose now that (4) holds. For every $1\leq i\leq n$, we consider $N_i=\prod_{j\neq i} P_j$, which is a characteristic subgroup of $G$ that is pro-$t$ embeddable into $P$ for all primes $t$. So $N_i$ is normal in $G$ and, by \cref{rec1}, $P/N_i\cong P_i$ is pro-$p_i$ embeddable into $G/N_i\cong P_i\rtimes Q$. By \cref{rec0},  the action of $Q$ on $P_i$ is trivial. We have proven this for all $i$, so the action of $Q$ on $P$ is trivial and hence we have (2).
\end{proof}

\begin{lem} \label{rec2} Let $G$ be a finite  group such that  for every subgroup $H \leq G$ and all primes $p$, the induced map $H_{\hp}\lrar G_{\hp}$ is injective. Then $G$ is nilpotent. 
\end{lem}
\begin{proof} We can prove this by induction on the order of $G$. When $|G|=1$ or, in fact, when it is a finite $p$-group for some $p$, the statement is trivial. Suppose then that $|G|$ is divisible by at least two primes $p$ and $q$. We consider a $p$-Sylow $P\leq G$ (which exists by the first Sylow theorem). Since $\{1\}$ is open in $P$, there must exist a normal subgroup $N\n G$ of $p$-power index such that $N\cap P=\{1\}$. Since $p$-subgroups of $G$ can be conjugated into $P$ (by the second Sylow theorem), then the order of $N$ is coprime with $p$. Consequently, $|P|=|G: N|$ and then $G=\lan N, P\ran$ is isomorphic to a semi-direct product $N\rtimes P$. We know that $|N|<|G|$ and then by the induction hypothesis, $N$ is nilpotent and so a direct product of its Sylow subgroups. By \cref{rec01}, we conclude that $G$ is nilpotent. 
\end{proof}

We have gathered all the ingredients  required to show \cref{reci}.
\begin{proof}[Proof of \cref{reci}] By \cref{rec1}, any finite quotient $Q$ of $G$ will still have the property that for every subgroup $H \leq Q$, the induced map $H_{\hp}\lrar Q_{\hp}$ is injective. So, by \cref{rec2}, every finite quotient of $G$ is nilpotent. Finally, by a result of Hirsch (see \cite[Theorem 2, Chapter 1]{Seg83}), a virtually polycyclic group has the property that all its finite quotients are nilpotent if and only if it is nilpotent itself. So $G$ is nilpotent. 
\end{proof}


\subsection{Pro-$p$ embeddability of polycyclic subgroups}
We are ready to show that in residually-(torsion-free nilpotent) groups, polycyclic groups are pro-$p$ embeddable (see \cref{pemb} below). Before doing so, we remind the reader of several classical properties of the class of virtually polycyclic groups (we refer the reader to the book of D. Segal \cite{Seg83} for a more thorough account). This class is closed under taking subgroups, quotients and extensions; they are noetherian (i.e. all their subgroups are finitely generated); and, crucially, that they are endowed with the useful notion of Hirsch length, which can be defined as follows. 

Given a virtually polycyclic group $G$ and a subnormal series $1=G_0\n G_1\n \cdots \n G_n=G$ such that the successive  quotients $G_{i+1}/G_i$ are cyclic, the number of subscripts $0\leq i\leq n-1$ such that  $G_{i+1}/G_i$ is infinite is an invariant of such a series (by the Jordan--H\" older theorem). We call this number the {\te Hirsch length} of $G$, denoted by $h(G)$. Recall that $h(G)=0$ if and only if $G$ is finite and that $h(G)=h(N)+h(G/N)$ for all normal subgroups $N\n G$. In particular, a surjective map between torsion-free virtually polycyclic groups $G_1\lrar G_2$ with $h(G_1)=h(G_2)$ must be an isomorphism. This will be important during the proof of the following.

\begin{thm}[\cref{RNsep}] \label{pemb} Let $G$ be a finitely generated residually-(torsion-free nilpotent) group and let $H\leq G$ be a virtually polycyclic subgroup. Then $H$ is nilpotent and, for all primes $p$, the induced map $H_{\hp}\lrar G_{\hp}$ is injective. 
\end{thm}
\begin{proof}

Let $\iota\po H\lrar G$ denote the injection and let $\pi^n$, for any positive integer $n$, denote the canonical projection $\pi^n\po G\lrar G/\ga_n G$. 

\begin{claim} \label{h0} There exists a positive integer $n$ such that the composition $f^n=\pi^n\circ \iota\po H\lrar G/\ga_n G$ is injective. 
\end{claim}
We will show by induction on the Hirsch length of $H$ the stronger statement that there exists a torsion-free nilpotent quotient $G/N$ such that the projection $\pi_N\po G\lrar G/N$ is injective on $H$. This is trivial if $h(H)=0$.  If $h(H)>0$, there exists a subgroup $H_1\leq H$  with $h(H_1)=h(H)-1$. By the induction hyphotesis, there exists a normal subgroup $N\n G$ such that the surjection  $\pi_N\po G\lrar G/N$ is injective on $H_1$. In particular, $h(\pi_N(H))\geq h(H_1)$. If $\pi_N$ is already injective on $H$ then we are done. Otherwise, there exists $1\neq x\in H\cap \ker \pi_N$.   Since $G$ is residually-(torsion-free nilpotent), there exists a normal subgroup $K\n G$ such that $\pi_K\po G\lrar G/K$ is injective on the subgroup $\lan x\ran\cong \Z$.  Now look at the induced surjection of torsion-free nilpotent groups $\pi_{N\cap K}(H)\lrar \pi_{N}(H)$. It is not injective because $\pi_{N\cap K}(x)\neq 1$ belongs to the kernel, so $h(\pi_{N\cap K}(H))\geq h(\pi_{N}(H))+1\geq h(H_1)+1=h(H)$. This implies that the surjection $\pi_{N\cap K}\po H\lrar \pi_{N\cap K}(H)$ of torsion-free virtually polycyclic groups is injective, completing the proof of \cref{h0}. In particular, $H$ is nilpotent.

Now we are ready to show that the induced map $H_{\hp}\lrar G_{\hp}$ is injective. Let $n\geq 1$ be an integer and consider the following natural commutative diagram. 
 \begin{equation} \label{magical}
     \begin{tikzcd}
       H \ar[dd, hookrightarrow, bend right=60, "f^n"'] \ar[d, hookrightarrow, "\iota"] \ar[r] & H_{\hp}  \ar[dd, bend left=60, "f^n_{\hp}"] \ar[d, "\iota_{\hp}"']\\
       G \ar[d, "\pi^n"] \ar[r]& G_{\hp} \ar[d, "\pi^n_{\hp}"'] \\
       G/\ga_n G \ar[r] & G_{\hp}/\ga_n G_{\hp}
     \end{tikzcd}
 \end{equation}
 
The horizontal arrows are the canonical pro-$p$ completion maps (the fact that the lowest arrow represents a pro-$p$ completion map follows from \cref{gakp}). We can choose an $n$ such that the map $f^n\po H\lrar G/\ga_n G$ is injective by \cref{h0}. Since $G/\ga_n G$ is finitely generated and nilpotent, it follows by \cref{pMalcev} that the induced $f_{\hp}$ is injective, too. In particular, $\iota_{\hp}$ must be injective, as we wanted.
\end{proof}

\subsection{Induced pro-$p$ topology on other infinite groups} \label{7.3} The pro-$p$ embeddability result of \cref{pemb} is a key component of the proof of \cref{Tits} on the 2-freeness of groups in the $p$-genus of a free or hyperbolic surface group. Here we introduce more context and motivation for the study of pro-$p$ separability properties and their applications, building on the setting of profinite topologies, which so far has been investigated more deeply.  There are geometric ways to produce families of LERF groups. Three big classes are limit groups \cite{Wil08}, locally quasi-convex hyperbolic virtually compact special groups (which includes   hyperbolic fundamental groups of graphs of free groups with cyclic edge subgroups \cite{Wis00}),  Seifert fibered spaces (either with boundary \cite[Theorem 5.1]{Hal49} or closed \cite[Theorem 4.1]{Sco78} \cite{Sco85}) and hyperbolic three manifolds (see \cite[Chapters 4--5]{Asc15} for complete referencing). 

Nevertheless, it seems very hard to use geometric methods to construct groups with good pro-$p$ embeddability or separability properties on its subgroups (as, for example, finitely generated nilpotent groups enjoy by \cref{pMalcev}). These difficulties appear when constructing parafree groups (for which, so far, we only have cyclic hierarchies \cite{And21}), when  classifying residually nilpotent 3-manifolds (see \cite[Question 2.12]{Wil08b}) or when it comes to find the right pro-$p$ analogue of limit groups \cite{Koc11}. On the other hand, there is research made on related  principles, such as the virtual $p$-efficiency of JSJ-decompositions of graph manifolds \cite{Asc13} and also of fibred 3-manifolds \cite[Theorem A]{Wilkes17}.

Apart from their own intrinsic interest, the motivation for studying these pro-$p$ separability principles is that they have a number of applications, for instance, on the residual properties and subgroup structure theory of pro-$p$ completions of ICE groups \cite{Koc11};  on the conjugacy $p$-separability of certain classes of three-manifold groups \cite[Theorems B and C]{Wilkes17}; on the profinite detection  of JSJ decompositions of graph manifolds \cite[Theorem A]{Wilkes18b}; and on the profinite rigidity of direct products of free and surface groups \cite[Theorem E]{Fru23}. These properties are useful and singular, and we consider that one should try to understand the following problem.

\begin{question} \label{introp2} Let $p$ be a prime. Are there finitely generated residually finite groups $G$ that contain non-abelian free subgroups and have the property of \cref{pMalcev}, i.e. that, for all finitely generated $H\leq G$, the pro-$p$ topology of $G$ induces on $H$ its full pro-$p$ topology?
\end{question}

\begin{remark}\label{plerf} Recall that a group is $p$-LERF if all its finitely generated subgroups are closed in the pro-$p$ topology. As we mentioned above, this property is too strong in the universe of infinite groups, since it does not even hold in $\Z$ (because $p\Z$ is not a closed subgroup in the pro-$q$ topology for any prime $q\neq p$). In fact, for similar reasons, a virtually polycyclic group is $p$-LERF if and only if it is a finite $p$-group. For a more general infinite group $G$ that is $p$-LERF, we can prove that all its finite quotients are $p$-groups (as in the proof of \cref{reci}). Equivalently, this means that its profinite completion is a pro-$p$ group. However, the property of being $p$-LERF is stronger than that of having a pro-$p$ group as profinite completion. In fact, Francoeur and Garrido \cite{Fra18} showed that all Sunic groups (with the exception of the infinite dihedral group) have the property that their profinite completion is a pro-2 group and, if they are not torsion, then they contain maximal finitely generated subgroups that are dense in the profinite topology (and hence they are not $2$-LERF).    Fern\'{a}ndez-Alcober, Garrido  and Uria-Albizuri \cite{Fer17} had already found examples with similar properties among GGS groups (standing for Grigorchuk--Gupta--Sidki groups). They provide examples of  finitely generated, residually finite and torsion-free groups $G$ whose profinite completion is a pro-$p$ group. Nevertheless, it is not clear to the author whether these groups are $p$-LERF.  
\end{remark}
The author is grateful to Alejandra Garrido for clarifying the picture that we describe in \cref{plerf}.

\section{The \texorpdfstring{$p$}{p}-genus of free or surface groups} \label{Titssec}

In this section we prove \cref{Titsintro}  from the introduction, although here he state it and prove it as part of a more general dichotomy on subgroups. 
 \begin{thm} \label{Tits} Let $p$ be a prime and let $S$ be a free or hyperbolic surface group. Suppose that $G$ is a finitely generated residually-$p$ group with $G_{\hp}\cong S_{\hp}$.  Assume that  $H\leq G$. The following holds.
 \begin{itemize}
     \item[(i)] If $H$ is amenable, then it is cyclic. 
     \item[(ii)] If $H$ is two-generated, then it is free. 
     \item[(iii)] If $H$ is finitely presented and non-cyclic, then it is large.
 \end{itemize}
 \end{thm}
 The property (i)  of this theorem will follow from (ii) by the following lemma. 
\begin{lemma} \label{cyclic} Let $G$ be a finitely generated  residually-(torsion-free polycyclic) group and let $H\leq G$ be locally-cyclic. Then $H$ is either trivial or infinite cyclic.
\end{lemma}
\begin{proof} If $H$ is finitely generated, then the conclusion follows immediately. We proceed by contradiction. Suppose, otherwise, that $H$ is not finitely generated. Take $1\neq h_1\in H$. Since $H$ is locally cyclic but not cyclic, there exists $a_2\notin \lan h_1\ran$ and also $h_2\in H$ such that $\lan h_2\ran=\lan h_1, a_2\ran$. So there exists $n_1\in \Z$  with $h_2^{n_1}=h_1$ and, due to the fact that $a_2\notin \lan h_1\ran$, we derive that $|n_1|>1$. Iterating this process, we show that then there must exist a sequence $(h_i)_{i\geq 1}$ in $H$ of non-trivial elements and a sequence of integers $(n_i)_{i\geq 1}$ with $|n_i|>1$ such that $h_{i+1}^{n_i}=h_i$ for all $i\geq 1$. Let $\phi\po G\lrar Q$ be a surjective map to a torsion-free polycyclic group $Q$ such that $\phi(h_1)\neq 1$ (and hence $\phi(h_i)\neq 1$  for all $i$). Now consider the ascending chain of subgroups
\[\lan \phi(h_1)\ran \sub \lan \phi(h_2)\ran \sub \cdots \sub \lan \phi(h_n)\ran \sub \cdots \]
of $Q$. Since $Q$ is noetherian, this sequence must stabilise. So there exists $k\geq 1$ such that $\phi(h_{k+1})\in \lan \phi(h_{k})\ran = \lan \phi(h_{k+1})^{n_k}\ran.$ Since $|n_k|>1$, this implies that $\phi(h_{k+1})$ has finite order in $Q$, and so $\phi(h_{k+1})=1$, which is a contradiction. 
\end{proof}

Before proving \cref{Tits}, we state the following well-known principle. 
\begin{lem} \label{2free} Let $S$ be a non-abelian free or surface group. Suppose that $\HH$ is a topologically two-generated closed subgroup of $S_{\hp}$. Then $\HH$ is a free pro-$p$ group.
\end{lem}
\begin{proof} If $S$ is free, this follows from \cref{1dim}. If $S$ is a hyperbolic surface group of genus $g$, then the abelianisation of $S$ has rank $2g\geq 4$. This implies that the image of $ \HH$ in the abelianisation of $S_{\hp}$ has infinite index and, in particular, $\HH$ already had infinite index in $S_{\hp}$. By \cref{Serre}, $\cd_p (\HH)\leq 1$. Lastly, again by \cref{1dim}, $\HH$ is a free pro-$p$ group.
\end{proof}

\begin{proof}[Proof of \cref{Tits}] We first explain how the property (i) follows from (ii). Remark that amenable groups cannot contain non-abelian free subgroups. So, if $H$ is amenable, then, by property (ii), $H$ is locally cyclic. In addition,  $G$ is a subgroup of $G_{\hp}$, which is residually-(torsion-free nilpotent) by \cref{tfnil2}. Finally, by \cref{cyclic}, $H$ is cyclic. 

Now we move on to prove (ii).  Let $S$ be a free or surface group and suppose that $G$ belongs to the $p$-genus of $S$. The statement is clear when $S\cong \Z$   because, in this case, $G\cong \Z$. Let us suppose that $S$ is a non-abelian free group or a hyperbolic surface group.  Let $H\leq G$ be a two-generated subgroup. We separately analyse two cases, depending on whether $H$ is abelian or not.
\begin{itemize}
    \item[(a)] Firstly, consider the case when $H$ is non-abelian. We consider $\ove{H}\sub G_{\hp}$. By  \cref{2free}, $\ove{H}$ is a two-generated free pro-$p$ group. Notice that $H\leq \ove{H}$ is not abelian, so the rank of $\ove{H}$ has to be exactly equal to two and then $\ove{H}$ is the free pro-$p$ group of rank 2. Let $x$ and $y$ be two elements that generate $H$, then $x$ and $y$ are topological generators of $\ove{H}$ and hence, by the Hopfian property, they are also free topological generators of $\ove{H}$. In particular, the group $H$ generated by $x$ and $y$ is free of rank two, as we wanted. 
    \item[(b)]  Lastly, suppose that $H$ is abelian. By \cref{pemb}, the induced map $H_{\hp}\lrar G_{\hp}$ is injective. Moreover, by \cref{2free},  $H_{\hp}$ is free. Hence $H_{\hp}$ is a cyclic pro-$p$ group and $H$ is cyclic.
\end{itemize}
 The previous shows (ii). Finally, we explain how (iii) follows. Let $H$ be a non-cyclic finitely presented subgroup. By (i), $H$ is not abelian, so $\ove H\leq S_{\hp}$ is not abelian either. The pro-$p$ group $\ove H$ must either be a non-abelian free pro-$p$ group or an open subgroup of $S_{\hp}$. In any case, $\ove H$ is isomorphic to the pro-$p$ completion of a non-abelian free or surface group $S'$. By \cref{bettigeq}, $\b(H)\geq \b(S')>0$. Thus $H$ is a finitely presented  residually-$p$ group with $\b(H)>0$ and, by \cite[Theorem 1.6]{Lac10}, it is large.  
\qedhere

\end{proof}

\begin{remark} \label{Titslimit} We should note that our methods can also establish parts (i) and (ii) of \cref{Tits} for other hyperbolic limit groups $S$. For example, one can take $S$ to be the free product $S=S_1*\cdots *S_n$, with each $S_i\cong F_i*_{u_i=\bar{u_i}} \bar{F_i}$ is the double of a free group $F_i$ along a relator ${u_i}\in F_i$ that satisfies the same assumptions as $r$ in \cref{tfnil2}. The pro-$p$ completion $S_{\hp}$ will be residually-(torsion-free nilpotent) by \cref{tfnil} and, by the proof of \cite[Theorem 7.3]{Koc11}, the 2-generated subgroups of $S_{\hp}$ are cyclic or free pro-$p$. So our proof of \cref{Tits} applies to such $S$ to show parts (i) and (ii).
\end{remark}

\section{The first \texorpdfstring{$L^2$}{L²}-Betti number as a pro-\texorpdfstring{$\CC$}{C} invariant} \label{Betti}
We will describe some scenarios in which the first $L^2$-Betti number of a group $G$, denoted by $\b(G)$, is a profinite invariant. This plays an important role in the proof of \cref{prods}. An standard reference about $L^2$-invariants is the book of L\" uck \cite{Luc02} and, in alignment with the approach taken in this section, the survey paper of Jaikin-Zapirain \cite{And19b}.

For torsion-free groups $G$ satisfying the Strong Atiyah conjecture,  Linnell  \cite{Lin93} shows that one can define $\b(G)$ as follows. Suppose that the group ring $\Q G$ of $G$ with rational coefficients has a   universal division ring of fractions $\Q G\inc \D_{\Q G}$. This provides a notion of dimension for all $\Q G$-modules  as follows: given a left $\Q G$-module $M$, we extend  scalars $\D_{\Q G} \otimes_{\Q G} M$ and take its linear-algebraic dimension as a left $ \D_{\Q G}$-module. We can define the first $L^2$-Betti number for these groups $G$ as 
\[\b(G)=\dim_{\D_{\Q G}} H_1(G; \D_{\Q G}).\]
Since locally indicable groups satisfy the Strong Atiyah conjecture by a result of Jaikin-Zapirain and L\' opez-\' Alvarez \cite{And192}, the previous is a definition of $\b$ for these groups. A general philosophy that explains why we work with  first $L^2$-Betti numbers instead of usual Betti numbers is that the former behave in a better way while containing similar information. For example, they are multiplicative in the following sense.

\begin{prop}[Theorem 1.35(9), \cite{Luc02}] \label{multi} Let $G$ be a finitely generated group and let $H$ be a finite-index subgroup. Then $ \b(H)=\b(G)\, |G: H|.$
\end{prop}

\subsection{Estimations using L\" uck approximation}
When it comes to relate the first $L^2$-Betti number $\b(G)$ of a group $G$ to its profinite completion $\hat G$, it is particularly helpful the characterisation of $\b(G)$ as a limit of (normalised) usual Betti numbers over a filtration of $G$ by normal finite-index subgroups.
A pioneering result in this direction is due to  L\" uck \cite{Luc94}.

\begin{thm}[L\" uck's approximation theorem] \label{Luck} Let  $G$ be a group of type $\FP_2$ and let $G=N_1>N_2>\dots >N_m>\dots $ be a sequence of finite-index normal subgroups  with $\bigcap_m N_m=1$. Then 
\[\lim_{m\rar \infty} \frac{b_1(N_m)}{|G:N_m|}= b_1^{(2)}(G).\]
\end{thm}
Importantly, for residually-finite groups that are not $\FP_2$, we still have a one-sided estimate of the same type. 

\begin{thm}[\cite{Luc11}] \label{Luckest} Let $G$ be finitely generated and let $G=N_1>N_2>\dots >N_m>\dots $ be a sequence of finite-index normal subgroups  with $\bigcap_m N_m=1$. Then
\[\limsup_{m\rar \infty} \frac{b_1(N_m)}{|G:N_m|}\leq  \b (G).\]
\end{thm}

We move on to the observation that groups $G$ satisfying the following inequality also have a L\" uck approximation theorem:
\begin{equation} \label{ine}
    b_1(G_0)\geq \b(G_0)+1\, \mbox{for every finite-index $G_0\leq G$.}
\end{equation}
Examples of groups $G$ with this property are limits of left-orderable amenable
groups in the space of marked group presentations (by a result of Osin \cite[Theorem 1.4]{Osi11}). 
Interestingly, the property of (\ref{ine}) can be proven in greater generality.

\begin{prop}[Jaikin-Zapirain \cite{And19}] \label{BettiAnd} Let $G$ be a finitely generated residually-(locally indicable and amenable) group. Then $\b(G)\leq b_1(G)-1$.
\end{prop} 
For our purposes, the interest behind \cref{BettiAnd} is that it allows to extend  the L\" uck approximation \cref{Luck} of $\b$  to other classes of groups. 

\begin{cor}\label{Luck2} Let $G$ be a finitely generated residually-(amenable and locally indicable) group. Let $G=N_1>N_2>\dots >N_m>\dots $ be a sequence of finite-index normal subgroups  with $\bigcap_m N_m=1$. Then
\begin{equation} \label{Luck2eq}
    \lim_{m\rar \infty} \frac{b_1(N_m)}{|G:N_m|}= \b (G).
\end{equation} 
\end{cor}
\begin{proof}
We know from \cref{Luckest} that
\[\limsup_{m\rar \infty} \frac{b_1(N_m)}{|G:N_m|}\leq  \b (G).\]
On the other hand, $b_1(N_m)\geq \b(N_m)+1=\b(G)|G: N_m|+1$ for all $m$ by \cref{BettiAnd} and the multiplicativity of $\b$ (\cref{multi}). So 
\[\liminf_{m\rar \infty} \frac{b_1(N_m)}{|G:N_m|}\geq  \b (G).\]
Hence the limit of (\ref{Luck2eq}) exists and equals $\b(G)$.
\end{proof}

\begin{defi} \label{Lu} We say that a finitely generated  residually finite group $G$ is $\Lu$ if for every descending sequence $G=N_1>N_2>\dots >N_m>\dots $ of finite-index normal subgroups $N_m$ with $\bigcap_m N_m=1$, it satisfies that
\[\lim_{m\rar \infty} \frac{b_1(N_m)}{|G:N_m|}= \b (G).\]
\end{defi}

We summarise this subsection by saying that both $\FP_2$ groups and finitely generated residually-(amenable and locally indicable) groups are $\Lu$ by \cref{Luck} and \cref{Luck2}.

\subsection{An invariant of the genus}

We are going to denote by $\CC$ a formation of finite groups contaning the formation of finite $p$-groups, denoted by $\CC_p$. If $G$ and $H$ are finitely generated groups with the same pro-$\CC$ completion then it is clear that $b_1(G)=b_1(H)$. This leads to the following corollary (which is a re-statement of \cite[Proposition 7.5]{Bri14}).

\begin{cor} \label{bettigeq}     Let $\CC$ be a formation of finite groups that includes $\CC_p$ for some $p$. Let $H$ be a $\Lu$ and residually-$\CC$ group. If $G$ admits a dense embedding into $H_{\hat \CC}$ then $\b(G)\geq \b (H)$. In particular, if $G$ is also $\Lu$, then $\b(G)=\b(H)$. 
\end{cor}

\section{Residually free and limit groups} \label{limitsec}

The purpose of this section is reviewing separability and structural properties of residually free groups  and a criteria to detect hyperbolic and limit groups among them. In \cref{special}, we will see to what extend some of these strong features are still true in some of the groups of the class $\Hia$. \cref{subprod} relates  the structure of finitely generated residually free groups to the structure of limit groups. 

\begin{prop}[\cite{Bau99}, Corollary 19; and \cite{Sel01}, Claim 7.5]\label{subprod} Let $G$ be a finitely generated residually free group. Then $G$ is a subdirect product of finitely many limit groups.
\end{prop}

A classical question attributed to Gromov asks whether one-ended hyperbolic groups must contain surface subgroups, and Wilton answers this question positively for many cyclic splittings in the following way. 

\begin{thm}[\cite{Wil18}, Theorem 6.1] \label{surfacesub} Let $G$ be the fundamental group of a graph of virtually free groups and virtually cyclic edge subgroups. Suppose that $G$ is hyperbolic and one-ended. Then $G$ contains a surface subgroup. 
\end{thm}

As proven in \cite[Corollary C]{Wil18}, the theorem above can be applied, in combination with Sela's hierarchy on limit groups \cite{Sel01}, to find surface subgroups in non-free limit groups. \cref{surfacesub} is slightly improved by Fruchter and the author \cite[Theorem A]{Fru23}.

\subsection{Criterion for the hyperbolicity of residually free groups}

We state a criterion that helps to recognise a hyperbolic group given a cyclic splitting. It is a consequence of the combination theorem for negatively curved groups of Bestvina--Feighn \cite{Bes92}.

\begin{thm}\label{comb} Let $(\G, Y)$ be a finite graph of groups with hyperbolic vertex groups and virtually cyclic edge subgroups. Then its fundamental group $\pi$ is hyperbolic if and only if it contains no Baumslag--Solitar subgroups. 
\end{thm}

This combination theorem was used by Sela \cite[Corollary 4.4. (iii)]{Sel01} to show that a limit group is hyperbolic if and only if it contains no $\Z^2$. This criterion can be extended to residually free groups as in the following result. 
 

\begin{prop}[\cite{Bri14a}, Lemma 4.15] \label{lim} Let $G$ be a finitely generated residually free group. Then $G$ is a hyperbolic limit group if and only if it contains no $\Z^2$.
\end{prop}

\subsection{Separability properties}

Wilton \cite{Wil08} proved that finitely generated subgroups of limit groups are quasiconvex and separable by establishing the following stronger property.

 \begin{defi} \label{retract} Let $A$ be a subgroup of $B$. We say that $A$ is a {\te virtual retract} of $B$ if there exists a finite-index subgroup $B_1$ of $B$ containing $A$ and a group homomorphism $r: B_1\lrar A$ such that $r$ restricts to the identity function on $A$. 
 \end{defi}

 \begin{thm}[\cite{Wil08}] Let $G$ be a limit group and let $H$ be a finitely generated subgroup. Then $H$ is a virtual retract of $G$. 
 \end{thm}
This principle was later extended to subgroups of type $\FP$ of finitely generated residually free groups by Bridson--Wilton \cite[Theorem B]{Bri08}. In the following section we treat other classes of groups with similar retraction properties (see \cref{VRspecial}). 


\section{Proof of \texorpdfstring{\cref{prods}}{Theorem E}} \label{prodsec}
In this final section we prove \cref{prods}. Let us review the statement.

\begin{thm}[Theorem F]  Let $G$ be a finitely generated residually free group and let $S$ be a free or surface group. Suppose that $\hat G\cong \hat S\times \hat \Z^n$. Then $G\cong S\times \Z^n$.
\end{thm}

We will state two lemmas that we will use during the proof of \cref{prods}. The first allows us to prove that the centre of a residually free group is a profinite invariant. 

\begin{lem}[\cite{Bau67a}, Lemma 4] \label{nocentre} Let $G$ be a residually free group. Then $G/Z(G)$ is a residually free group with trivial centre. 
\end{lem}

The second lemma is a profinite criterion that detects when short exact sequences split. It was used by Wilton--Zalesskii \cite{Wilton17} to distinguish circle bundles over the surface by looking at finite quotients of their fundamental group (later on used in dimension 4 by Jiming--Zixi \cite[Proposition 25]{Ma22}).
\begin{lem}[\cite{Wilton17}, Lemma 8.3]\label{splitt} Let $n\geq 0$ be an integer and let $S$ be a surface group. Let
\begin{equation} \label{split}
     \begin{tikzcd}
       1 \ar[r]  & \Z^n \ar[r]  & G \ar[r] & S \ar[r] & 1
       \end{tikzcd}
 \end{equation}
 be a central extension of $S$ by the group $\Z^n$. Suppose that the induced short sequence in profinite completions
 \begin{equation} \label{split0}
     \begin{tikzcd}
       1 \ar[r]  & \hat \Z^n \ar[r]  & \hat G \ar[r] & \hat S \ar[r] & 1
       \end{tikzcd}
 \end{equation}
  is also exact  and that it splits. Then the sequence (\ref{split}) splits, too. 
\end{lem}

 We remark that the exactness of (\ref{split0}) is always ensured if (\ref{split}) is exact, by the goodness of $S$ (as observed by Serre \cite[Chapter 1, Section 2, Exercise 2(b)]{Ser97}).  However, we stated this as an assumption for clearness, as it is also enough for our application. Regarding this lemma of Serre, we also refer the reader to \cite[Proposition 2.4]{Lor08} for a precise proof and a converse. 

\begin{proof}[Proof of \cref{prods}] Let $S$ be a free or surface group and let $G$ be a finitely generated residually free group in the finite genus of $S\times \Z^n$. If $S$ is abelian, then it is clear that $G\cong S\times  \Z^n$. So we can suppose that $S$ is hyperbolic and non-abelian. In particular, $\b(S)>0$ and $\hat S$ has trivial centre. We denote by $ \iota\po Z(G) \inc G$ the canonical injection. 

\begin{claim} \label{cl1} The centre $Z(G)$ is a finitely generated subgroup of $G$ and the induced map $\hat \iota\po \hat{Z(G)}\lrar \hat G$ is  injective.
\end{claim}
To prove the claim, first recall that $Z(\hat G)= \hat \Z^n$ is a retract of $\hat G$. Denoting this retraction by $r\po \hat G\lrar Z(\hat G)$,  we have the following natural commutative diagram:

   \begin{equation} \label{magic2}
     \begin{tikzcd}
       Z(G) \ar[dd, rightarrow, bend right=60, "f"'] \ar[d, hookrightarrow, "\iota"] \ar[r, hookrightarrow] &  Z( \hat G) \ar[dd, bend left=60, "\id"] \ar[d, hookrightarrow]\\
       G \ar[d, "p"] \ar[r, hookrightarrow]& \hat G \ar[d, "r"'] \\
       G/[G, G] \ar[r, "g"] & Z(\hat G),
     \end{tikzcd}
 \end{equation}
for some naturally defined $f$ and $g$ from the canonical maps $\iota, p$ and $r$. We can read off from the diagram that $g\circ f$ is injective, so $f$ is injective, too. This implies that $Z(G)$ is a subgroup of the finitely generated abelian group $G/[G, G]$, so it is finitely generated, too. Now we consider a second natural commutative diagram:
   \begin{equation} \label{magic3}
     \begin{tikzcd}
       Z(G) \ar[dd, hookrightarrow, bend right=60, "f"'] \ar[d, hookrightarrow, "\iota"] \ar[r, hookrightarrow] &  \hat{Z(G)} \ar[dd, bend left=60, "\hat f"] \ar[d, "\hat \iota "']\\
       G \ar[d, "p"] \ar[r, hookrightarrow]& \hat G \ar[d, "\hat p"'] \\
       G/[G, G] \ar[r, hookrightarrow] & \widehat{G/[G, G]}.
     \end{tikzcd}
 \end{equation}
We know that $f$ is injective, so $\hat f$ is also injective (since finitely generated abelian groups are LERF). Thus $\hat \iota $ is injective and \cref{cl1} is proven. 

\begin{claim} \label{cl12} The profinite group $ Z(\hat G)/\overline{Z(G)}$ is the profinite completion of a finitely generated abelian group $A$.
\end{claim}

The image of $\hat \iota $ can be alternatively described as the closure of $Z(G)\inc \hat G$ and it is easy to see that it is a normal subgroup of $\hat G$ that is contained in $Z(\hat G)=\hat\Z^n$. We consider the group $H=G/Z(G)$, which is a finitely generated residually free group with trivial centre by \cref{nocentre}. By  \cref{propquo}, $\hat H$ is naturally isomorphic to $\hat G/\overline{Z(G)}\cong  \hat S\times \left( Z(\hat G)/\overline{Z(G)}\right)$.
So $Z(\hat G)/\overline{Z(G)}$  appears as a direct factor of  $\hat H/[\hat H, \hat H]\cong\left(\hat S/[\hat S, \hat S]\right)\times \left( Z(\hat G)/\overline{Z(G)}\right)$, where  $\hat H/[\hat H, \hat H] $  and $\hat S/[\hat S, \hat S]$ are the profinite completions of $H/[H, H]$ and $S/[S, S]$, respectively. This shows \cref{cl12}.

\begin{claim} \label{cl2} The image of the map $\hat \iota\po \hat{Z(G)}\lrar \hat G$ is $Z(\hat G)$. 
\end{claim} 

We denote by $p_1$ and $p_2$ the projections of $\hat H$ onto its first and second direct factor according to the decomposition $ \hat H\cong \hat S \times Z(\hat G)/\overline{Z(G)}$. \cref{cl2} follows if the second factor is trivial. Let us suppose that it is not trivial to reach a contradiction. We study two separate cases. 
\begin{itemize}
    \item[(i)] Suppose that $Z(\hat G)/\overline{Z(G)}$ is infinite. By \cref{cl12}, there exists an infinite finitely generated abelian group $A$ such that $Z(\hat G)/\overline{Z(G)}\cong \hat A$. Hence $H$ and $K=S\times A$ have the same profinite completion. Observe that $K$ is $\Lu$ and that $\b(K)=0$. Since $H$ is residually free, it is also $\Lu$ by \cref{Luck2}, so it follows from \cref{bettigeq} that $\b(H)= \b(K)=0$. This implies that the restriction of $p_1$ to $H$ cannot be injective (otherwise, we would see $H$ as a dense subgroup of $\hat S$ and it would follow from \cref{bettigeq} that $\b(H)\geq \b(S)>0$, contradicting the fact that $\b(H)=0$). So $H$ intersects $\ker p_1=Z(\hat G)/\overline{Z(G)}=Z(\hat H)$ non-trivially and then $Z(H)\neq 1$, which contradicts \cref{nocentre}. 
    \item[(ii)] Suppose that $Z(\hat G)/\overline{Z(G)}$ is finite. Then we can still prove that $H$ intersects non-trivially the direct factor $Z(\hat G)/\overline{Z(G)}$, leading to the same contradiction as before. To verify this, let us denote $H_1=H\cap \hat S$. This is a finite-index subgroup of $H$ with $\hat H_1\cong \hat S$ by \cref{bij}. So $H_1$ is also $\Lu$ and $\b(H_1)=\b (S)$ by \cref{bettigeq}. Observe that $p_1(H_1)$ is isomorphic to $H_1$, that this is finite-index in $p_1(H)$, and that both $p_1(H_1)$ and $p_1(H)$ are dense in  $\hat S$. Let us denote the index $k=|p_1(H): p_1(H_1)|$. Using the multiplicativity of $\b$ (\cref{multi}) and the estimation of \cref{bettigeq}, it follows that \[\b(p_1(H))\geq \b(S)=\b(H_1)=\b(p_1(H_1))=k\cdot \b(p_1(H)).\] Since $\b(S)$ is positive, we derive from the previous inequalities that $k=1$ and hence $p_1(H)=p_1(H_1)$. From this, it is immediate to see that $Z(\hat G)/\overline{Z(G)}=p_2(H)=p_2(H\cap Z(\hat G)/\overline{Z(G)})$, and then $H\cap Z(\hat G)/\overline{Z(G)}\neq 1$, as we wanted. 
\end{itemize}
This proves \cref{cl2}. Since $H=G/Z(G)$ is a finitely generated residually free group and $\hat H\cong \hat S$ by \cref{cl2}, it follows from the work of Wilton \cites{Wil18, Wil21}  that $H\cong S$. 
\begin{claim} The short exact sequence  
 \begin{equation} \label{magic4}
     \begin{tikzcd}
       1 \ar[r]  & Z(G) \ar[r]  & G \ar[r] & G/Z(G) \ar[r] & 1
       \end{tikzcd}
 \end{equation}
 splits. In particular, $G\cong Z(G)\times G/Z(G)\cong \Z^n\times S$.
\end{claim}
The proof of the claim would complete the proof of \cref{prods}. This claim is obvious when $S$ is free and, when $S$ is a surface group, it is a direct consequence of \cref{splitt}.
\end{proof}

\section{Cohomological dimension and Poincar\'e duality} \label{cohsec}

Here we  recall the finiteness properties of groups that we shall use. 
We  denote by $R$ a unital and associative ring (not necessarily commutative). 

\begin{defi} \label{finiteness} We say that an $R$-module $M$ is of type  $\FP_n(R)$ if there exists a projective resolution 
\[\lrar P_{n+1}\lrar P_n\lrar \cdots \lrar P_1\lrar P_0\lrar M\lrar 0\]
with $P_i$ finitely generated for $0\leq i\leq n$. Moreover, we say that an $R$-module $M$ is of type $\FP_{\infty}(R)$ if it is of type $\FP_n(R)$ for all $n$. Analogously, we say that a group $G$ is of type $\FP_n(R)$ (resp. of type $\FP_{\infty}(R)$) if the trivial $R G$-module
$R$ is of type $\FP_n(R)$ (resp. of type $\FP_{n}(R)$ for all $n$).
\end{defi}

\subsection{Profinite groups} We refer the reader to \cite[Chapter 2, Section 7]{Neu08} for details about the cohomology theory of profinite groups with profinite coefficient modules. Computations on profinite group cohomology can be reduced to analogous computations on abstract group cohomology via Serre's  fundamental notion of goodness \cite[Section I. 2. 6]{Ser97}. A group is said to be {\te cohomologically good} (or simply {\te good}) if for every finite discrete $\hat G$-module $A$, the natural homomorphism 
\[H^k(\hat{G}; A)\lrar H^k(G; A),\]
induced by $G\lrar \hat{G}$, is an isomorphism. 
A group $G$ is said to be {\te subgroup separable} or {\te LERF} (standing for {\te locally extended residually finite}) if every finitely generated subgroup $H\leq G$ is closed in the profinite topology. This implies that  the topology induced on $H$ from the  profinite topology of $G$ is the full profinite topology of $H$. Hence, by \cref{embed}, if $G$ is LERF, then the natural map $\hat{H}\lrar \bar{H}\sub \hat{G}$ is an isomorphism of profinite groups for all finitely generated $H\leq G$. 
  Given a profinite group $\GG$, the \textit{completed group algebra  of $\GG$ over $\Z_p$} is a profinite ring defined as the inverse limit of the usual group rings $\left(\Z/p^i\Z\right)[\GG/\UU]$,
 where $i\geq 0$ and $\UU$ ranges over open normal subgroup of $\GG$. We say that  $\GG$ is of type $p$-$\FP_n$  if the $\llbracket\hat \Z \GG\rrbracket$-module $\hat \Z$ (resp. the $\llbracket\Z_p\GG\rrbracket$-module $\Z_p$) is of type $\FP_n$.
 
\begin{prop}[\cite{And18}, Proposition 3.1] If $\Ga$ is an abstract $\fpi$ and good group then its profinite completion $\hat{\Ga}$ is of type $p$-$\fpi$ for every prime $p$. 
\end{prop}
Recall we say a module $M$ is {\te simple} if and only if it has exactly two submodules.
\begin{defi} \label{cdp} We say that a profinite group $\GG$ has $p$-cohomological dimension $n=\cd_p(\GG)$ if $n$ is the largest non-negative integer $m$ such that  $H^{m+1}(\GG, A)=0$ for all simple discrete $\llbracket\hat{\Z}\GG\rrbracket$-modules that are annihilated by $p$.
\end{defi}
Recall from \cite[Corollary 1 of Chapter 3, and Proposition 21]{Ser97} that if $\GG_p$ is the $p$-Sylow of a profinite group $\GG$, then $\cd_p(\GG)=\cd_p(\GG_p)$. Good groups $G$ enjoy the convenient property that
\[\cd_p(\hat{G})=\cd_p (G)\, \, \,\mbox{for all primes $p$.}\] In particular, a free profinite group has $\cd_p= 1$ and profinite surfaces have $\cd_p= 2$. We should remark that focusing on $p$-primary modules in the \cref{cdp} of the ``cohomological dimension function'' $\cd_p$ allows to avoid the pathology of $H^2(\hat \Z; \Z)$ being non-zero (as this is isomorphic to $H^1(\hat \Z, \Q/\Z)\cong \Q/\Z$  \cite[Chapter 1, Section 3.2]{Ser97}).
The dimension function $\cd_p$ enjoys nice ``geometric'' properties. For example, by  Shapiro's lemma given a closed subgroup $\HH$ of a profinite group $\GG$, we have $\cd_p(\HH)\leq \cd_p(\GG)$ (see \cite[Theorem 7.3.1]{Rib00}). Furthermore, $\cd_p$ is sometimes additive with respect to short exact sequences (as reflected by \cite[Theorem 1.1]{Wei04}). 

A classical characteristic feature of groups of finite cohomological dimension is the following well-known lemma.

\begin{lem} Let $\GG$ be a profinite group with $\cd_p (\GG)<\infty$ for all primes $p$. Then $\GG$ is torsion-free. 
\end{lem}
This is esentially the only tool we have to recognise torsion-freeness of an abstract group from the profinite completion. Interestingly, Lubotzky shows that being torsion-free is not a profinite property (\cite[Proposition 1.5]{Lub14}). The following is the analogous theorem of Stallings' \cite{Sta68} and Swan's  \cite{Swa69} for the category of pro-$p$ groups.
\begin{thm}[Theorem 7.7.4, \cite{Rib00}] \label{1dim} Let $\GG$ be a pro-$p$ group with $\cd_p(\GG)\leq 1$. Then $\GG$ is a free pro-$p$ group.
\end{thm}


\subsection{Poincar\'e Duality and the recognition of surfaces} 
 Eckman and M\" uller \cite{Eck80}  proved that surface groups are the only  Poincar\'e duality groups of dimension two. Before discussing the analogous principles for profinite groups, we recall the following definition.
 
\begin{defi} Let $\GG$ be a profinite group of type $p$-$\fpi$. We say that $\GG$ is a {\te Poincar\'e duality group at $p$} of dimension $n$ (written as $\PD^n$ at $p$) if $\cd_p(\GG)=n$ and 
\begin{align*}
    &H^i(\GG, \llbracket\Z_p \GG\rrbracket)=0 \, \,\, \,  \mbox{for $i\neq n$}, \\
    &H^i(\GG, \llbracket\Z_p \GG\rrbracket)\cong \Z_p \, \,\, \,  \mbox{(as abelian groups).}
\end{align*}
These denote the {\it continuous cochain cohomology} groups in the sense of \cite[Definition 2.7.1]{Neu08}.
\end{defi}
 
 A pro-$p$ group $\GG$ is $\PD_1$  if and only if $\GG\cong \Z_p$ (see \cite[Example 4.4.4]{Sym00}) and Poincar\'e duality pro-$p$ groups of dimension 2 are known to coincide with the class of {\te Demushkin} groups (see \cite[Section I.4.5, Example 2]{Ser97}), which is a class of one-relator pro-$p$ groups that was  classified by Demushkin, Serre and Labute \cites{Dem61, Dem63, Ser62, Lab67}. However, this question  is poorly understood for general profinite groups. Similarly, it is not known whether a finitely generated residually finite group whose profinite completion is $\PD_2$ at a prime $p$ is itself $\PD_2$. Wilton \cite{Wil21} showed that if such group is additionally residually free then it is a surface group, which we restate in \cref{surfacedet}. This has been recently generalised in \cite{Jai23} with a different argument. 
 In order to restate Wilton's criterion, we require the following profinite variant of Strebel's theorem \cite{Str77} on infinite-index subgroups of $\PD_n$ abstract groups (which can be found in \cite[Page 44, exercise 5(b)]{Ser97} and \cite[Chapter 3, Section 7, Exercise 3]{Neu08}). 

\begin{lem}[Serre] \label{Serre} Let $\GG$ be a profinite $\PD_n$ group at a prime $p$ and let $\HH$ be a closed subgroup such that $p^{\infty}$ divides the index $|\GG: \HH|$. Then $\cd_p(\HH)<n$. 
\end{lem}

The previous lemma was of fundamental importance in the work of Wilton--Zalesskii \cite{Wil19} on the profinite detection of prime and JSJ decompositions of three-manifolds. 

\begin{remark} When $\GG$ is a pro-$p$ group, \cref{Serre} only requires $\HH$ to have infinite index. Nevertheless, for general profinite groups $\GG$, the conclusion of the theorem would not be true by only requiring infinite index. In fact, from \cref{Zhat}, we know that $\hat \Z$, which is $\PD_1$ at $p$, contains closed subgroups isomorphic to $\hat \Z$ of index $\prod_p p$.
\end{remark}

The most natural source of examples of profinite Poincar\'e duality groups consists of profinite completions of abstract Poincar\'e duality groups, as the following result from \cite[Theorem 4.1]{Koc08} ensures. 
\begin{thm} \label{PDcom} Let $\Ga$ be an abstract cohomologically good $\PD_n$ group. Then $\hat{\Ga}$ is a profinite $\PD_n$ group at every prime and $\Ga_{\hp}$ is $\PD_n$ at $p$.
\end{thm}

Finally, we are ready to give a more refined criterion to detect surfaces, used by Wilton \cite{Wil21} to establish \cref{Rem2} when $G$ is a limit group.

\begin{prop}\label{surfacedet} Let $p$ be a prime and let $G$ be a torsion-free residually finite group with separable cyclic subgroups such that its profinite completion $\hat G$ is $\PD_2$ at $p$. Suppose that $G$ contains a surface subgroup as a virtual retract. Then $G$ is a surface group.
\end{prop}

\begin{proof}  By assumption, there exists a finite-index subgroup $L$ of $G$ containing a surface group $\pi_1\Sigma$ as a retract. Let $K$ be the kernel of the retraction $L\lrar \pi_1\Sigma$. Our aim is to prove that $K=1$. This way, we would have that $L\cong \pi_1\Sigma$ and  that the  torsion-free group $G$ would be  virtually a surface group. So the conclusion would follow from Kerckhoff's solution to the Nielsen realisation problem for surfaces. 

Hence, it remains to show that $K=1$. We proceed by contradiction. Suppose that $K\neq 1$. The retraction $\pi_1 \Sigma\inc L\lrar \pi_1\Sigma$ induces a retraction $\hat{\pi_1\Sigma}\inc \hat{L}\lrar \hat{\pi_1\Sigma}$ at the level of the profinite completions. Take $1\neq k\in K$. Since $G$ is torsion-free and cyclic subgroups are separable, then $\lan k\ran\cong \Z$ and the supernatural number $p^{\infty}$ divides the order of $\overline{\lan k\ran} \cong \hat \Z$. Since $k$ is contained in the kernel $\KK$ of the retraction $\hat{L}\lrar \hat{\pi_1 \Sigma}$, then $|\KK|=|\overline{\lan k\ran} ||\KK: \overline{\lan k\ran} |$ is also divisible by $p^{\infty}$. We now claim that $p^{\infty}$ also divides the index $|\hat{L}:\hat{\pi_1\Sg}|$. This way, it would also divide the index $|\hat{G}:\hat{\pi_1\Sg}|$ and, by \cref{Serre}, the $p$-cohomological dimension of $\hat{\pi_1\Sg}$ should be at most $\cd_p(\hat{G})-1=1$, which is a contradiction. 

Thus, it simply remains to prove that $p^{\infty}$ divides $|\hat{L}\!:\!\hat{\pi_1\Sg}|$. From the comments above, we know that, for every $n\geq 1$, there exists a characteristic open subgroup $\UU\no \KK$ such that $p^n$ divides $|\KK\!:\!\UU|$. We observe that $\HH=\UU \cdot \hat{\pi_1\Sigma}$ is an open subgroup of $\hat L$. After taking a normal open subgroup $\HH_0\no \hat L$ contained in $\HH$ (for example, its core), we are going to check that $p^n$ divides $|\hat L: \HH_0\cdot  \hat{\pi_1\Sg}|$. For this, is suffices to notice that the previous index is divisible by
\[|\hat L: \HH\cdot \hat{\pi_1\Sg}|=|\KK \cdot \hat{\pi_1\Sigma}: \UU \cdot \hat{\pi_1\Sigma}|=|\KK: \UU|,\]
since $\UU\no \KK$ and $\KK\cap \hat{\pi_1\Sigma}=\{1\}$. We chose $\UU$ so $p^n$ divides $|\KK: \UU|$, so this proves that $p^{\infty}$ divides $|\hat{L}: \hat{\pi_1\Sg}|$, as we wanted.
\end{proof}

For the purpose of proving \cref{Remthm} using \cref{surfacedet}, we remark that the profinite completion of a surface $\hat {\pi_1\Sg}$ is $\PD_2$ at $p$ for all primes $p$.

\section{Hyperbolic and special groups}\label{special}

The study of special cube complexes was  initiated by Haglund--Wise in their seminal paper \cite{Hag08} and plays a fundamental role in the solution 
 of many results about three-manifolds, such as Agol's proof of the virtual Haken conjecture \cite{Ago13}. We refer the reader to the book of Aschenbrenner--Friedl--Wilton \cite{Asc15} for more examples and discussion. A group is said to be {\te special} (resp. {\te compact special}) if it is the fundamental group of a non-positively curved (resp. compact) cube complex that satisfies Haglund--Wise's special condition on its hyperplanes. We will not delve into this definition because we are only interested in the retraction properties of these groups (as reflected by \cref{VRspecial}) and their consequences in the proof of \cref{Remthm}. Nevertheless, in order to offer a more precise idea of what these groups look like, we recall the following alternative description. 

\begin{thm}[\cite{Hag08}, Theorem 4.2]  A group $G$ is special (resp. compact special) if and only if it is a subgroup (resp. a quasiconvex subgroup) of a right-angled Artin group. 
\end{thm}

As anticipated in the introduction,  many groups from $\Hia$ will be hyperbolic and virtually compact special.  

\begin{prop} \label{hvcs} Let $G$ be a group belonging to $\Hia$ with the following properties: 
\begin{enumerate}
    \item It is torsion-free.
    \item It contains no Baumslag--Solitar subgroups.
    \item Its abelian subgroups are cyclic. 
\end{enumerate}
Then $G$ is  hyperbolic and virtually compact special.
\end{prop}
\begin{proof} The group $G$ is hyperbolic by \cref{comb}. Now the conclusion that $G$ is virtually compact special follows directly from Wise's classification of   hyperbolic groups with a malnormal
quasi-convex hierarchy \cite[Theorem 11.2]{Wis21}.
\end{proof}

 From the proof  of \cite[Theorem 7.3]{Hag08}, we know that quasi-convex subgroups of hyperbolic compact special groups are virtual retracts (recall \cref{retract}).

\begin{thm} \label{VRspecial} Let $G$ be a hyperbolic virtually compact special group and let $H\leq G$ be a quasi-convex subgroup. Then there exists a finite-index subgroup $H' $ of $H$ such that $H'$ is a virtual retract of $G$. 
\end{thm}

We can use Wilton's partial solution to Gromov's surface conjecture (\cref{surfacesub}) to find surface subgroups in certain groups of $\Hia$ that will be relevant for the proof of \cref{Remthm}.

\begin{prop} \label{surfacesub2} Let $G$ be as in \cref{hvcs}. Suppose that $G$ is not free. Then $G$ contains a surface subgroup as a virtual retract. 
 \end{prop}
 
 \begin{proof} By \cref{hvcs}, $G$ is hyperbolic and virtually compact special. So, by \cref{VRspecial}, in order to conclude that $G$ contains a hyperbolic surface subgroup as a virtual retract, it is enough to ensure that $G$ contains a quasi-convex surface subgroup. Furthermore, we know that for finitely generated subgroups of a hyperbolic group, being quasi-convex and quasi-isometrically embedded is equivalent. Hence, it suffices to ensure that there exists a quasi-isometrically embedded surface subgroup in $G$ to prove \cref{surfacesub2}, which we prove by induction on the level of the hierarchy.  

  Let us start with the level 0, that is, when $G$ is a finitely generated residually free group. Since $G$ contains no $\Z^2$, then $G$ is limit by \cref{lim} and the statement follows directly from \cref{surfacesub}. Now we establish the inductive step by distinguishing two cases.
\begin{itemize}
    \item  Suppose that $G=A *_C B$ (resp. $G=A _{C, \theta}$ for some injection $\theta: C\lrar A$), where $A$ and $B$ (resp. $A$) are groups in $\Hia$ for which the statement holds and $C$ is cyclic. Since both $A$ and $B$ are quasi-isometrically embedded in $G$ (for example, by \cite[Theorem 1.2]{Kap01}), then $G$ will contain a quasi-isometrically embedded surface if at least one  of $A$ or $B$ does. So it remains to ensure the conclusion in the case when $A$ and $B$ (resp. $A$) are free. Viewing $G$ as the fundamental group of a graph of two free groups (resp. one free group) with infinite cyclic edge groups, we apply a variant of  Shenitzer's theorem \cite[Theorem 18]{Wil12} to deduce that, if $G$ is not one-ended, then some vertex will split freely relative to its incident edge groups. In particular, applying this result a finite number of times (and taking out an infinite cyclic free factor of $G$ each time), we deduce that $G$ is isomorphic to the free product $F* G'$ for a free group $F$ and a one-ended group $G'=A' *_C B'$ (resp. $G=A' _{C, \theta}$ for the injection $\theta: C\lrar A'$), where $A'$ and $B'$ are free factors of $A$ and $B$ (resp. $A'$ is a free factor of $A$). Since $G$ was hyperbolic, and $G'$ is a quasi-convex subgroup, then $G'$ is also hyperbolic. So $G'$ is a hyperbolic one-ended fundamental group of a graph of free groups with cyclic edges. By \cref{surfacesub}, the group $G'$ (and hence also $G$) contains a quasi-convex surface subgroup. 
    \item Suppose that $G$ has a finite-index subgroup $H\in \Hia$ for which the statement is already true. Since $G$ is torsion-free and not free, then $H$ is not free by \cite[Theorem 2]{Sta68}. So $H$  (and hence also $G$) must contain a surface subgroup that is a virtual retract. \qedhere
\end{itemize}
 \end{proof}

\section{Pro-\texorpdfstring{$p$}{p} detection of algebraic and topological properties}\label{detection}
 We recall that by \cref{And}, groups in the finite, soluble or nilpotent genus of  a free or surface group $S$, do belong to the $p$-genus of $S$ for all primes $p$. So all the results here about any $p$-genus of $S$ include all these genera.
\subsection{L\" ucky groups}
It is generally not well understood if being L\" ucky (\cref{Lu}) is a profinite invariant (see \cref{Luinv} below). However, Jaikin-Zapirain \cite{And212} can answer affirmatively to this question for the genera of free and surface groups as a consequence of the following. 

\begin{prop} \label{Luck3} Let $S$ be a free or surface group and let $G$ be an abstract subgroup of $S_{\hp}$. Then $G$ is $\Lu$. 
\end{prop}
\begin{proof} By \cref{tfnil} and \cref{tfnil2}, $G$ is residually-(torsion-free nilpotent). By \cref{Luck2}, $G$ is $\Lu$.
\end{proof}

It is not difficult to see that the previous $G$ have to be RFRS (there is, in fact, a stronger statement in \cite[Proposition 4.1]{And212}). 

\begin{question} \label{Luinv} Let $\Ga$ and $\La$ be finitely generated residually finite groups with $\hat \Ga\cong \hat \La$. If $\Ga$ is $\Lu$, does it follow that $\La$ is $\Lu$?
\end{question} 
This relates the question of whether being RFRS is a profinite property. 
\begin{question} Let $\Ga$ and $\La$ be finitely generated residually finite groups with $\hat \Ga\cong \hat \La$. If $\Ga$ is RFRS, does it follow that $\La$ is RFRS?
\end{question}

\subsection{Hyperbolic and special groups}
Here we see how \cref{Tits} allows us to study when groups from $\Hia$ (\cref{defhia}) or one-relator groups can be proven to be hyperbolic and virtually special by looking at their pro-$p$ completion.

 \begin{cor} \label{specialdet} Let $p$ be a prime and let $G$ be a group in the class $\Hia$ that belongs to the $p$-genus of a free or hyperbolic surface group. Then $G$ is hyperbolic and virtually compact special. 
 \end{cor}
 \begin{proof} By   \cref{Tits}, $G$ is torsion-free, contains no Baumslag--Solitar groups and contains no non-cyclic abelian subgroups. Hence, by  \cref{hvcs}, $G$ is hyperbolic virtually compact special.
 \end{proof}

Recently, Linton \cite[Theorem 8.2]{Lin22} established that all one-relator groups with no negative immersions are hyperbolic and virtually special. By \cite[Theorem 1.3]{Lou22}, a one-relator group that is two-free has no negative immersions. Hence  the following  is a consequence of these results and \cref{Tits}.

\begin{cor}[\cref{hvspecialintro}] \label{onerel} Let $G$ be a one-relator group that belongs to the  $p$-genus of a free group or a hyperbolic surface group. Then $G$ is hyperbolic and virtually special.
\end{cor}


\section{On the profinite rigidity of free and surface groups} \label{freesurface}

We are now ready to establish the profinite rigidity of free and surface groups within the family $\Hia$. The proof of \cref{Remthm} is divided in \cref{Remfree} and \cref{Remsur}.

\begin{thm} \label{Remfree} Let $G$ be a group belonging to the class $\Hia$ and let $F$ be a finitely generated free group. If $G$ is residually finite and $\hat G\cong \hat{F}$, then $G\cong F$. 
\end{thm}

\begin{proof} Since free groups are distinguished from each other by their abelianisation, it suffices to show that $G$ is free. Let $G$ be a group in $\Hia$ that belongs to the finite genus of a free group $F$. By \cite[Theorem 1.1]{And212}, $G$ also belongs to the $p$-genus of $F$. Consequently, by \cref{specialdet}, $G$ is hyperbolic virtually compact special. To prove that $G$ is free, we proceed by contradiction. Suppose that is $G$ is not free.  Then, by \cref{surfacesub2}, there exists a surface subgroup $\pi_1 S$  that it is a virtual retract of $G$. In particular, $\hat {\pi_1 S}\inc \hat G\cong \hat F$. This implies by the goodness of the surface $\pi_1 S$ and the monotonicity of $\cd_p$ the following contradiction:
\[2=\cd_p(\pi_1 S)=\cd_p (\hat {\pi_1 S})\leq \cd_p( \hat F)=1. \qedhere\]
\end{proof}

\begin{thm} \label{Remsur} Let $G$ be a group belonging to the class $\Hia$ and let $\pi_1\Sg$ be the fundamental group of a closed orientable surface. If $G$ is residually finite and $\hat G\cong \hat{\pi_1\Sg}$, then $G\cong \pi_1 \Sg$. 
\end{thm}

\begin{proof} Since surface groups are distinguished from each other by their abelianisation, it suffices to show that $G$ is a surface group. Let $G$ be a group in $\Hia$ that belongs to the finite genus of a free group $\pi_1\Sg$. By \cite[Theorem 1.1]{And212}, $G$ also belongs to the $p$-genus of $\pi_1\Sg$. Consequently, by \cref{specialdet}, $G$ is hyperbolic virtually compact special. Since $G$ cannot be free, it follows from \cref{surfacesub2} that there exists a surface subgroup $\pi_1 S$ that it is a virtual retract of $G$. Lastly, by \cref{surfacedet}, $G$ is a surface group. 
\end{proof}


\textsc{Mathematical Institute, University of Oxford, Radcliffe Observatory, Andrew Wiles Building, Woodstock Rd, Oxford OX2 6GG} \\
\textit{E-mail address:} \href{mailto:morales@maths.ox.ac.uk}{\texttt{morales@maths.ox.ac.uk}}

\bibliography{biblio.bib}

\end{document}